\documentclass[reqno,11pt]{amsart}
\usepackage{amssymb}

\usepackage{amsmath}
\usepackage{times}
\usepackage[dvips]{graphicx}

\newtheorem{theorem}{Theorem}

\newtheorem{corollary}[theorem]{Corollary}

\newtheorem{lemma}[theorem]{Lemma}

\newtheorem{proposition}[theorem]{Proposition}

\begin{document}
\title{A concentration inequality and a local law for the sum of two random
matrices}
\author{Vladislav Kargin}
\thanks{Department of Mathematics, Stanford University, CA 94305;
kargin@stanford.edu}
\date{June 2011}
\maketitle

\begin{center}
\textbf{Abstract}
\end{center}

\begin{quotation}
Let $H_{N}=A_{N}+U_{N}B_{N}U_{N}^{\ast }$ where $A_{N}$ and $B_{N}$ are two $%
N$-by-$N$ Hermitian matrices and $U_{N}$ is a Haar-distributed random
unitary matrix, and let $\mu _{H_{N}},$ $\mu _{A_{N}}$, $\mu _{B_{N}}$ be
empirical measures of eigenvalues of matrices $H_{N}$, $A_{N}$, and $B_{N}$,
respectively. Then, it is known (see \cite{pastur_vasilchuk00}) that for
large $N$, the measure $\mu _{H_{N}}$ is close to the free convolution of
measures $\mu _{A_{N}}$ and $\mu _{B_{N}}$, where the free convolution is a
non-linear operation on probability measures. The large deviations of the
cumulative distribution function of $\mu _{H_{N}}$ from its expectation have
been studied by Chatterjee in \cite{chatterjee07}. In this paper we improve
Chatterjee's concentration inequality and show that it holds with the rate
which is quadratic in $N.$ \newline
In addition, we prove a local law for eigenvalues of $H_{N_{N}},$ by showing
that the normalized number of eigenvalues in an interval approaches the
density of the free convolution of $\mu _{A}$ and $\mu _{B}$ provided that
the interval has width $\left( \log N\right) ^{-1/2}.$
\end{quotation}

\section{Introduction}

If $A$ and $B$ are two Hermitian matrices with a known spectrum, it is a
classical problem to determine all possibilities for the spectrum of $A+B.$
The problem goes back at least to H. Weyl (\cite{weyl12}). Later, Horn (\cite%
{horn62}) suggested a list of inequalities which must be satisfied by
eigenvalues of $A+B$, and recently, Knutson and Tao (\cite{knutson_tao99})
using earlier ideas by Klyachko, proved that this list is complete.

For large matrices, it is natural to consider the probabilistic analogue of
this problem, when matrices $A$ and $B$ are ``in general position''. Namely,
let $H_{N}=A_{N}+U_{N}B_{N}U_{N}^{\ast },$ where $A_{N}$ and $B_{N}$ are two
fixed $N$-by-$N$ Hermitian matrices, and $U_{N}$ is a random unitary matrix
with the Haar distribution on the unitary group $\mathcal{U}\left( N\right)
. $ Then, the eigenvalues of $H_{N}$ are random and we are interested in
their joint distribution.

Let $\lambda _{1}^{\left( A\right) }\geq \ldots \geq \lambda _{N}^{\left(
A\right) }$ denote eigenvalues of $A_{N},$ and define the \emph{spectral
measure }of $A_{N}$ as $\mu _{A_{N}}:=N^{-1}\sum_{k=1}^{N}\delta _{\lambda
_{k}^{\left( A\right) }}.$ Define $\mu _{B_{N}}$ and $\mu _{H_{N}}$
similarly, and note that $\mu _{H_{N}}$ is random even if $\mu _{A_{N}}$ and 
$\mu _{B_{N}}$ are non-random. What can be said about relationship of $\mu
_{A_{N}},$ $\mu _{B_{N}},$ and $\mu _{H_{N}}$?

An especially interesting case occurs when $N$ is large. This case was
investigated by Voiculescu (\cite{voiculescu91}) and Speicher (\cite%
{speicher93}) who found that as $N$ grows $\mu _{H_{N}}$ approaches $\mu
_{A_{N}}\boxplus \mu _{B_{N}},$ where $\boxplus $ denotes \emph{free
convolution}, a non-linear operation on probability measures introduced by
Voiculescu in his studies of operator algebras. Their proofs are based on
calculating traces of large powers of matrices and use ingenious
combinatorics. Later, Pastur and Vasilchuk (\cite{pastur_vasilchuk00})
applied the method of Stieltjes transforms to this problem and extended the
results of Speicher and Voiculescu to measures with unbounded support.

It appears natural to ask the question about deviations of $\mu _{H_{N}}$
from $\mu _{A_{N}}\boxplus \mu _{B_{N}}.$

In order to illuminate the issues that arise, suppose first that we place $N$
points independently on a fixed interval $\left[ a,b\right] ,$ each
according to a measure $\nu .$ Let the number of points in a sub-interval $I$
be denoted $\mathcal{N}_{I}$. Then, $\mathcal{N}_{I}$ is a sum of
independent Bernoulli variables and satisfies the familiar central limit law
and large deviation estimates. In particular, 
\begin{equation}
\Pr \left\{ \left| \frac{\mathcal{N}_{I}}{N\left| I\right| }-\mathbb{E}%
\left( \frac{\mathcal{N}_{I}}{N\left| I\right| }\right) \right| >\delta
\right\} \sim c_{1}\exp \left[ -c_{2}\delta ^{2}N\right]
\label{large_deviations_nonmatrix}
\end{equation}%
for large $N.$

A remarkable fact is that for random points corresponding to eigenvalues of
classical random matrix ensembles, the asymptotic is different and given by
the formula%
\begin{equation}
\Pr \left\{ \left| \frac{\mathcal{N}_{I}}{N\left| I\right| }-\mathbb{E}%
\left( \frac{\mathcal{N}_{I}}{N\left| I\right| }\right) \right| >\delta
\right\} \sim c_{1}\exp \left[ -c_{2}f\left( \delta \right) N^{2}\right] .
\label{large_deviations_matrix}
\end{equation}%
Intuitively, there is a repulsion force between eigenvalues which makes
large deviations of $\mathcal{N}_{I}$ much more unlikely for large $N$.

For classical ensembles this fact was rigorously shown in a more general
form in \cite{ben_arous_guionnet97}. Later, this result was extended to
matrices of the form $A_{N}+sX_{N},$ where $A_{N}$ is an Hermitian $N$-by-$N$
matrix and $X$ is an Hermitian Gaussian $N$-by-$N$ matrix; see for an
explanation Sections 4.3 and 4.4 in \cite{anderson_guionnet_zeitouni10}.

The fluctuations of eigenvalues of matrices $H_{N}=A_{N}+U_{N}B_{N}U_{N}^{%
\ast }$ were considered by Chatterjee in \cite{chatterjee07}. By an
ingenious application of the Stein method he proved that for every $x\in 
\mathbb{R}$, 
\begin{equation*}
\Pr \left\{ \left| \mathcal{F}_{H_{N}}\left( x\right) -\mathbb{E}\mathcal{F}%
_{H_{N}}\left( x\right) \right| >\delta \right\} \leq 2\exp \left[ -c\delta
^{2}\frac{N}{\log N}\right] ,
\end{equation*}%
where $\mathcal{F}_{H_{N}}\left( x\right) :=N^{-1}\mathcal{N}_{\left(
-\infty ,x\right] }$ denotes the cumulative distribution function for
eigenvalues of $H_{N},$ symbol $\mathbb{E}$ denotes the expectation with
respect to the Haar measure, and $c$ is a numeric constant. Note that the
rate in this estimate is sublinear in $N,$ hence the estimate is weaker than
(\ref{large_deviations_matrix}). In fact, it is even weaker than the
estimate in (\ref{large_deviations_nonmatrix}) because of the logarithmic
factor $(\log N)^{-1}$, and therefore it does not contain any evidence of
the repulsion between eigenvalues.

The first main result of this paper is an improvement of this estimate and
is as follows.

\textbf{Assumption }$A1.$ The measure $\mu _{A_{N}}\boxplus \mu _{B_{N}}$ is
absolutely continuous everywhere on $\mathbb{R},$ and its density is bounded
by a constant $T_{N}.$

\begin{theorem}
\label{theorem_main}Suppose that Assumption $A1$ holds. Let $\mathcal{F}%
_{H_{N}}$ and $\mathcal{F}_{\boxplus ,N}$ be cumulative distribution
functions for the eigenvalues of $H_{N}=A_{N}+U_{N}B_{N}U_{N}^{\ast }$ and
for $\mu _{A_{N}}\boxplus \mu _{B_{N}}$, respectively. Then, for all $N\geq
\exp \left( \left( c_{1}/\delta \right) ^{4/\varepsilon }\right) ,$%
\begin{equation}
P\left\{ \sup_{x}\left| \mathcal{F}_{H_{N}}\left( x\right) -\mathcal{F}%
_{\boxplus ,N}\left( x\right) \right| >\delta \right\} \leq \exp \left[
-c_{2}\delta ^{2}N^{2}\left( \log N\right) ^{-\varepsilon }\right] ,
\label{inequality_main}
\end{equation}%
where $c_{1},c_{2}$ are positive and depend only on $K_{N}:=\max \left\{
\left\| A_{N}\right\| ,\left\| B_{N}\right\| \right\} $, $T_{N}$, and $%
\varepsilon \in \left( 0,2\right] $.
\end{theorem}

Up to a logarithmic factor, the rate in this inequality is proportional to $%
N^{2},$ which is consistent with the possibility that the eigenvalues of
matrix $H_{N}=A_{N}+U_{N}B_{N}U_{N}^{\ast }$ repulse each other.

With respect to Assumption $A1,$ it is pertinent to note that if $\mu
_{A_{N}}\left( \left\{ x\right\} \right) <1/2$ and $\mu _{B_{N}}\left(
\left\{ x\right\} \right) <1/2$ for every $x\in \mathbb{R}$ (i.e., if the
multiplicity of every eigenvalue of $A_{N}$ and $B_{N}$ is less than $N/2$),
then $\mu _{A_{N}}\boxplus \mu _{B_{N}}$ has no atoms (see Theorem 7.4 in %
\cite{bercovici_voiculescu98}). Moreover, since $\mu _{A_{N}}$ and $\mu
_{B_{N}}$ are atomic, the results of \cite{belinschi08} imply that the
density of $\mu _{A_{N}}\boxplus \mu _{B_{N}}$ is analytic (i.e., in $%
C^{\infty }$ class) everywhere on $\mathbb{R}$ where it is positive. In
particular, Assumption $A1$ holds.

If Assumption $A1$ is relaxed, then it is still possible to prove a result
similar to the result in Theorem \ref{theorem_main}. Namely, if $\mathcal{%
\mu }_{\boxplus ,N}$ is absolutely-continuous at the endpoints of interval $%
I,$ then it is possible to show that for all sufficiently large $N,$ 
\begin{equation}
P\left\{ \left| \frac{\mathcal{N}_{I}}{N\left| I\right| }-\mathcal{\mu }%
_{\boxplus ,N}\left( I\right) \right| >\delta \right\} \leq \exp \left[
-c_{2}\delta ^{2}N^{2}\left( \log N\right) ^{-\varepsilon }\right] .
\label{inequality_alternative}
\end{equation}%
\ Indeed, the only place where Assumption $A1$ is used is when the distance
between $\mathcal{F}_{H_{N}}$ and $\mathcal{F}_{\boxplus ,N}$ is estimated
in terms of the distance between $m_{H}(z)$ and $m_{\boxplus ,N}(z)$ and
this is done by using Bai's theorem. In order to prove (\ref%
{inequality_alternative}), the original proof should be modified by using
techniques from the proof of Corollary 4.2 in \cite{erdos_schlein_yau09}
instead of Bai's theorem. In this paper, however, we choose to concentrate
on the proof of inequality (\ref{inequality_main}).

In addition, if Assumption $A1$ fails and $x$ is an atom of $\mu
_{A_{N}}\boxplus \mu _{B_{N}}$ then by Thm 7.4 in \cite%
{bercovici_voiculescu98} there exist $x_{A}$ and $x_{B}$ such that $%
x_{A}+x_{B}=x,$ and 
\begin{equation*}
\mu _{A_{N}}\left( \left\{ x_{A}\right\} \right) +\mu _{B_{N}}\left( \left\{
x_{B}\right\} \right) -1=\mu _{A_{N}}\boxplus \mu _{B_{N}}\left( \left\{
x\right\} \right) .
\end{equation*}%
These $x_{A}$ and $x_{B}$ are eigenvalues of $A_{N}$ and $%
U_{N}B_{N}U_{N}^{\ast }$ with multiplicities $\mu _{A_{N}}\left( \left\{
x_{A}\right\} \right) N$ and $\mu _{B_{N}}\left( \left\{ x_{B}\right\}
\right) N$, respectively. Hence, by counting dimensions and using the fact
that eigenspaces of $A_{N}$ and $U_{N}B_{N}U_{N}^{\ast }$ are in general
position, we conclude that with probability $1,$ $x_{A}+x_{B}$ is an
eigenvalue of $H_{N}$ with multiplicity 
\begin{equation*}
\left( \mu _{A_{N}}\left( \left\{ x_{A}\right\} \right) +\mu _{B_{N}}\left(
\left\{ x_{B}\right\} \right) -1\right) N.
\end{equation*}%
Hence, if $x$ is an atom of $\mu _{A_{N}}\boxplus \mu _{B_{N}}$, then we
have the exact equality 
\begin{equation*}
\mu _{H_{N}}\left( \left\{ x\right\} \right) =\mu _{A_{N}}\boxplus \mu
_{B_{N}}\left( \left\{ x\right\} \right) .
\end{equation*}

These considerations suggest that perhaps Assumption $A1$ can be eliminated
or weakened as a condition of Theorem \ref{theorem_main}$.$

Our main tools in the proof of Theorem \ref{theorem_main} are the Stieltjes
transform method and standard concentration inequalities applied to
functions on the unitary group.

In the first step, we establish the $N^{2}$ rate for large deviations of the
Stieltjes transform of $\mu _{H_{N}},$ which we denote $m_{H_{N}}\left(
z\right) .$ This follows from results in \cite{anderson_guionnet_zeitouni10}
and the fact that the Stieltjes transform of $\mu _{H_{N}}$ is Lipschitz as
a function of $U_{N}$ and its Lipschitz constant can be explicitly estimated.

It is not possible to prove a concentration inequality for $\mathcal{F}%
_{H_{N}}\left( x\right) $ by a similar method because for some $x$ this
function is not Lipschitz in $U_{N}.$ An alternative is to use an inequality
by Bai (Theorem \ref{theorem_bai} in this paper), which gives a bound on $%
\sup_{x}\left| \mathcal{F}_{H_{N}}\left( x\right) -\mathbb{E}\mathcal{F}%
_{H_{N}}\left( x\right) \right| $ in terms of $\sup_{x}\left|
m_{H_{N}}\left( z\right) -\mathbb{E}m_{H_{N}}\left( z\right) \right| ,$
where $z=x+i\eta .$ However, the second term in this inequality depends on
smoothness of $\mathbb{E}\mathcal{F}_{H_{N}}\left( x\right) $, which is
difficult to establish.

Instead, we show that $\sup_{x}\left| \mathbb{E}m_{H_{N}}\left( z\right)
-m_{\boxplus ,N}\left( z\right) \right| $ is small for $\eta :=\mathrm{Im}%
z>c/\sqrt{\log N}.$ (Here $m_{\boxplus ,N}\left( z\right) $ denote the
Stieltjes transform of $\mu _{A_{N}}\boxplus \mu _{B_{N}}$.) This estimate
allows us to use Bai's inequality and estimate $\sup_{x}\left| \mathcal{F}%
_{H_{N}}\left( x\right) -\mathbb{E}\mathcal{F}_{\boxplus ,N}\left( x\right)
\right| $ in terms of the sum of $\sup_{x}\left| m_{H_{N}}\left( z\right) -%
\mathbb{E}m_{H_{N}}\left( z\right) \right| $ and $\sup_{x}\left| \mathbb{E}%
m_{H_{N}}\left( z\right) -m_{\boxplus ,N}\left( z\right) \right| ,$ which
are both small. The benefit of this change is that smoothness of $\mathcal{F}%
_{\boxplus ,N}\left( x\right) $ is easier to establish than the smoothness
of $\mathbb{E}\mathcal{F}_{H_{N}}\left( x\right) .$ In our case it is
guaranteed by Assumption $A1.$

For large $\mathrm{Im}z,$ the difference $|\mathbb{E}m_{H_{N}}(z)-m_{%
\boxplus ,N}(z)|$ can be estimated by applying Newton's iteration method (as
perfected by Kantorovich in \cite{kantorovich48}) to the Pastur-Vasilchuk
system for $\mathbb{E}m_{H_{N}}(z)$. Namely, we use $m_{\boxplus ,N}(z)$ as
the starting point for this method and show that for sufficiently large $N$
the difference of the solution of the system, $\mathbb{E}m_{H_{N}}\left(
z\right) ,$ and the starting point is less than any fixed $\delta >0$.

This method fails for small $\mathrm{Im}z.$ We use a modification of
Hadamard's three circle theorem (\cite{hardy15}) in order to estimate the
difference $\left| \mathbb{E}m_{H_{N}}\left( z\right) -m_{\boxplus ,N}\left(
z\right) \right| $ in the region close to the real axis.

Theorem \ref{theorem_main} implies the following local law result. Let $%
N_{\eta }\left( E\right) $ denote the number of eigenvalues of $H_{N}$ in an
interval of width $2\eta $ centered at $E,$ and let $\varrho _{\boxplus
,N}\left( E\right) $ denote the density of $\mu _{A_{N}}\boxplus \mu
_{B_{N}} $ at $E.$

\begin{theorem}
\label{theorem_local_law} Suppose that $\eta =\eta \left( N\right) $ and $1/%
\sqrt{\log N}\ll \eta \ll 1$. Let assumption $A1$ hold with $T_{N}=T$.
Assume also that $\max \left\{ \left\| A_{N}\right\| ,\left\| B_{N}\right\|
\right\} \leq K$ for all $N.$ Then, for all sufficiently large $N,$ 
\begin{equation*}
P\left\{ \sup_{E}\left| \frac{\mathcal{N}_{\eta }\left( E\right) }{2N\eta }%
-\varrho _{\boxplus ,N}\left( E\right) \right| \geq \delta \right\} \leq
\exp \left( -c\delta ^{2}\frac{\left( \eta N\right) ^{2}}{\left( \log
N\right) ^{2}}\right) ,
\end{equation*}%
where $c>0$ depends only on $K$ and $T$.
\end{theorem}

(Here the notation $( N) \ll g( N) $ means that $\lim_{N\rightarrow \infty
}g\left( N\right) /f(N)=+\infty .$)

The plan of the rest of the paper is as follows. We start in Section \ref%
{section_notation} by establishing our notation. Section \ref%
{Section_concentration_Stieljtes} provides a large deviation estimate for
the Stieltjes transform of $\mu _{H_{N}}$ and a related function. In Section %
\ref{Section_PV_System}, we use this estimate to bound error terms in the
Pastur-Vasilchuk system, which we re-derive for reader's convenience.
Section \ref{Section_stability_PV_system} is devoted to estimating $\left| 
\mathbb{E}m_{H_{N}}\left( z\right) -m_{\boxplus ,N}\left( z\right) \right| $
in the region where $\mathrm{Im}z\geq \eta _{0}$, and Section \ref%
{section_hardy_theorem} is concerned with estimating it in the region $%
\mathrm{Im}z\gg 1/\sqrt{\log N}$. Section \ref{section_proof_main_theorem}
completes the proof of our two main theorems. Several concluding remarks are
made in Section~\ref{section_conclusion}.

\section{\textbf{Definitions and} \textbf{Notations}}

\label{section_notation}

We define $H_{N}=A_{N}+U_{N}B_{N}U_{N}^{\ast }.$ The \emph{spectral measure}
of $H_{N}$ is $\mu _{H_{N}}:=N^{-1}\sum_{k=1}^{N}\delta _{\lambda
_{k}^{(H)}},$ where $\lambda _{k}^{(H)}$ are eigenvalues of $H,$ counted
with multiplicity. Its \emph{cumulative distribution function} is denoted $%
\mathcal{F}_{H_{N}}\left( x\right) :=\mu _{H_{N}}\left( \left( -\infty ,x%
\right] \right) .$ The number of eigenvalues of $H_{N}$ in interval $I$ is
denoted $\mathcal{N}_{I}:=N\mu _{H_{N}}\left( I\right) ,$ and $\mathcal{N}%
_{\eta }\left( E\right) :=\mathcal{N}_{\left( E-\eta ,E+\eta \right] }$
denotes the number of eigenvalues in the interval of width $2\eta $ centered
at $E.$

The \emph{resolvent} of $H_{N}$ is defined as $G_{H}\left( z\right) :=\left(
H_{N}-z\right) ^{-1}.$ Similarly, $G_{A}\left( z\right) :=\left(
A_{N}-z\right) ^{-1}$ and $G_{B}\left( z\right) :=\left( B_{N}-z\right)
^{-1}.$ (For brevity, we will omit the subscript $N$ in the notation for
resolvents and Stieltjes transforms.)

The \emph{Stieltjes transform} of $H_{N}$ is defined as 
\begin{equation*}
m_{H}\left( z\right) :=N^{-1}\mathrm{Tr\,}G_{H}\left( z\right) =\int_{%
\mathbb{R}}\frac{\mu _{H_{N}}\left( d\lambda \right) }{\lambda -z},
\end{equation*}%
where $\mathrm{Tr}$ denotes the usual matrix trace$.$ The Stieltjes
transforms of $A_{N}$ and $B_{N}$ are defined similarly, e.g., $m_{A}\left(
z\right) =N^{-1}\mathrm{Tr\,}G_{A}\left( z\right) .$ More generally, if $\mu 
$ is a probability measure, then its Stieltjes transform is defined as 
\begin{equation*}
m_{\mu }\left( z\right) :=\int_{\mathbb{R}}\frac{\mu \left( d\lambda \right) 
}{\lambda -z}.
\end{equation*}

In addition, we define the following quantities: 
\begin{equation*}
f_{B}\left( z\right) :=N^{-1}\mathrm{Tr}\left( U_{N}B_{N}U_{N}^{\ast }\frac{1%
}{H_{N}-z}\right)
\end{equation*}%
and$.$%
\begin{equation*}
f_{A}\left( z\right) :=N^{-1}\mathrm{Tr}\left( A_{N}\frac{1}{H_{N}-z}\right)
\end{equation*}

Next, we define the free convolution. Consider the following system: 
\begin{eqnarray}
m\left( z\right) &=&m_{A}\left( z-S_{B}\left( z\right) \right) ,
\label{free_convolution_system} \\
m\left( z\right) &=&m_{B}\left( z-S_{A}\left( z\right) \right) ,  \notag \\
z+\frac{1}{m\left( z\right) } &=&S_{A}\left( z\right) +S_{B}\left( z\right) ,
\notag
\end{eqnarray}%
where $m\left( z\right)$, $S_{A}\left( z\right) ,$ $S_{B}\left( z\right) $
are unknown functions.

\begin{proposition}
\label{prop_uniqueness}There exists a unique triple of analytic functions $%
m(z),S_{A}(z),S_{B}(z)$ that are defined in $\mathbb{C}^{+}=\left\{ z:%
\mathrm{Im}z>0\right\} ,$ satisfy system (\ref{free_convolution_system}),
and have the following asymptotics as $z\rightarrow \infty $: 
\begin{eqnarray}
m\left( z\right) &=&-z^{-1}+O\left( z^{-2}\right) ,
\label{asymptotic_conditions} \\
S_{A,B}\left( z\right) &=&O\left( 1\right) .  \notag
\end{eqnarray}%
Moreover, the function $m\left( z\right) $ maps $\mathbb{C}^{+}$ to $\mathbb{%
C}^{+}$ and the functions $S_{A,B}\left( z\right) $ map $\mathbb{C}^{+}$ to $%
\mathbb{C}^{-}=\left\{ z:\mathrm{Im}z<0\right\} .$
\end{proposition}

Prop. \ref{prop_uniqueness} implies that the first function in this triple, $%
m_{\boxplus ,N}(z)$, is the Stieltjes transform of a probability measure.
This measure is called the \emph{free convolution} of measures $\mu _{A_{N}}$
and $\mu _{B_{N}}$ and denoted $\mu _{A_{N}}\boxplus \mu _{B_{N}}.$ (For
shortness, we will sometimes write this measure as $\mu _{\boxplus ,N}$.)
The two other functions in this triple, $S_{A}\left( z\right) $ and $%
S_{B}\left( z\right) ,$ are called \emph{subordination functions}.

\textbf{Proof of Prop. \ref{prop_uniqueness}:} The uniqueness of the
solution of system (\ref{free_convolution_system}) was proved in Prop. 3.3
in (\cite{pastur_vasilchuk00}). However, it appears that their proof does
not show that the solution exists everywhere in the upper half-plane. We
prove the existence and uniqueness differently, by establishing a one-to-one
correspondence between solutions of (\ref{free_convolution_system}) and
certain objects in free probability theory. After this correspondence is
established, the existence, uniqueness and claimed properties of the
solution follow from the corresponding properties of the free probability
objects.

Recall that in the traditional definition of free convolution (see \cite%
{voiculescu_dykema_nica92}), one defines the $R$-transform of measure $\mu
_{A}$ by the formula $R_{A}\left( t\right) =m_{A}^{\left( -1\right) }\left(
-t\right) -1/t,$ where $m_{A}^{\left( -1\right) }$ is the functional inverse
of $m_{A},$ chosen in such a fashion that $R_{A}\left( t\right) $ is
analytic at $t=0.$ The function $R_{B}\left( t\right) $ is defined
similarly. Then, one proves that $R=R_{A}+R_{B}$ is the $R$-transform of a
probability measure, and one calls this measure the free convolution of $\mu
_{A}$ and $\mu _{B}.$ In fact, this definition of free convolution is
equivalent to the definition we have given above.

Indeed, let $m_{\boxplus ,N}$ be the Stieltjes transform of $\mu
_{A_{N}}\boxplus \mu _{B_{N}}$ as it is usually defined, that is, let it
equal the functional inverse of $R+1/t$ multiplied by $-1.$ By definition of 
$R_{A}$, the first equation of (\ref{free_convolution_system}) can be
written equivalently as%
\begin{equation*}
S_{B}\left( z\right) =z+\frac{1}{m_{\boxplus ,N}\left( z\right) }%
-R_{A}\left( -m_{\boxplus ,N}\left( z\right) \right) ,
\end{equation*}%
which we can use as a definition of $S_{B}\left( z\right) .$ This definition
holds only for sufficiently large $z.$ However, by the results of Biane (%
\cite{biane98b}), $S_{B}\left( z\right) $ can be analytically continued to
the whole of $\mathbb{C}^{+}.$ If we write the second equation in a similar
form, add them together, and use the equality $R=R_{A}+R_{B},$ then we get: 
\begin{eqnarray*}
S_{A}\left( z\right) +S_{B}\left( z\right) &=&2z+\frac{2}{m_{\boxplus
,N}\left( z\right) }-R\left( -m_{\boxplus ,N}\left( z\right) \right) \\
&=&z+\frac{1}{m_{\boxplus ,N}\left( z\right) },
\end{eqnarray*}%
which is the third equation of system (\ref{free_convolution_system}). By
analytic continuation it holds everywhere in $\mathbb{C}^{+}.$ This shows
that if $m_{\boxplus ,N}\left( z\right) ,$ $S_{A}\left( z\right) ,$ and $%
S_{B}\left( z\right) $ are defined using the traditional definition of free
convolution, then they satisfy system (\ref{free_convolution_system}). In
particular this shows the existence of the solution of (\ref%
{free_convolution_system}) as a triple of analytic functions defined
everywhere in $\mathbb{C}^{+}.$

Conversely, if $m_{\boxplus ,N}\left( z\right) ,$ $S_{A}\left( z\right) ,$
and $S_{B}\left( z\right) $ satisfy (\ref{free_convolution_system}) with
asymptotic conditions (\ref{asymptotic_conditions}), then in a neighborhood
of infinity we can write 
\begin{eqnarray*}
R_{A}\left( -m_{\boxplus ,N}\left( z\right) \right) &=&m_{A}^{\left(
-1\right) }\left( m_{\boxplus ,N}\left( z\right) \right) +1/m_{\boxplus
,N}\left( z\right) \\
&=&z-S_{B}\left( z\right) +1/m_{\boxplus ,N}\left( z\right) ,
\end{eqnarray*}%
where the first line is the definition of $R_{A}$ and the second uses the
first equation of (\ref{free_convolution_system}). If we write a similar
expression for $R_{B}\left( -m_{\boxplus ,N}\left( z\right) \right) ,$ add
them together, and use the third equation of (\ref{free_convolution_system}%
), then we find that 
\begin{equation*}
R(-m_{\boxplus ,N}(z))=z+1/m_{\boxplus ,N}(z).
\end{equation*}%
This shows that $m_{\boxplus ,N}\left( z\right) $ satisfies the same
functional equation as the Stieltjes transform of the free convolution
measure defined in the traditional fashion. Since their power expansions at
infinity are the same, these functions coincide. In particular, this shows
that the solution of (\ref{free_convolution_system}) is unique as a triple
of analytic functions in $\mathbb{C}^{+}$ that satisfy asymptotic conditions
(\ref{asymptotic_conditions}).

Finally, the claimed properties of $m_{\boxplus ,N}\left( z\right) $ and $%
S_{A,B}\left( z\right) $ follow from the properties of the Stieltjes
transform of a probability measure and of the subordination functions. The
latter were established by Biane in (\cite{biane98b}). $\square $

We denote the cumulative distribution function of $\mu _{A_{N}}\boxplus \mu
_{B_{N}}$ as $\mathcal{F}_{\boxplus ,N}$ and its density (when it exists) as 
$\varrho _{\boxplus ,N}.$

The integration over $U$ using the Haar measure will be denoted as $\mathbb{E%
}.$ (This operation is often denoted as $\left\langle \cdot \right\rangle $
in the literature.) Correspondingly, $P\left( \omega \right) $ denotes the
Haar measure of event $\omega .$

We will usually write $z=E+i\eta ,$ where $E$ and $\eta $ denote the real
and imaginary parts of $z.$ We will also use the following notation: 
\begin{equation*}
\Omega _{\eta _{0},c}=\left\{ z\in \mathbb{C}:\mathrm{Im}z\geq \eta _{0},%
\mathrm{Im}z\geq c\mathrm{Re}z\right\} .
\end{equation*}

\section{Concentration for the Stieltjes transform and associated functions}

\label{Section_concentration_Stieljtes}

The main result of this section is the following large deviation estimates
for $m_{H}\left( z\right) $ and $f_{B}\left( z\right) .$

\begin{proposition}
\label{theorem_m_large_deviations}Let $z=E+i\eta $ where $\eta >0.$ Then,
for a numeric $c>0$ and every $\delta >0,$ 
\begin{equation}
P\left\{ \left| m_{H}\left( z\right) -\mathbb{E}m_{H}\left( z\right) \right|
>\delta \right\} \leq \exp \left( -\frac{c\delta ^{2}\eta ^{4}}{\left\|
B\right\| ^{2}}N^{2}\right) ,  \label{concentration_for_m}
\end{equation}%
and%
\begin{equation}
P\left\{ \left| f_{B}\left( z\right) -\mathbb{E}f_{B}\left( z\right) \right|
>\delta \right\} \leq \exp \left[ -\frac{c\delta ^{2}\eta ^{4}}{\left\|
B\right\| ^{4}}N^{2}/\left( 1+\frac{\eta }{\left\| B\right\| }\right) ^{2}%
\right] .  \label{concentration_for_f}
\end{equation}
\end{proposition}

\textbf{Proof: }The first claim of this proposition follows directly from
Corollary 4.4.30 in \cite{anderson_guionnet_zeitouni10}. The second claim
can be obtained by a modification of the proof of this Corollary. For the
convenience of the reader we give a short proof of both claims.

Both claims are consequences of the Gromov-Milman results about the
concentration of Lipschitz functions on Riemannian manifolds (\cite%
{gromov_milman83}). In a small neighborhood of identity matrix, all unitary
matrices can be written as $U=e^{iX}$, where $X$ is Hermitian. We identify
the space of Hermitian matrices $X$ with $T\mathcal{U},$ the tangent space
to $\mathcal{U}\left( N\right) $ at point $U.$ By left translations this
identification can be extended to the tangent space at any point of $%
\mathcal{U}\left( N\right) .$ Define an inner product norm in $T\mathcal{U}$
by the formula $\left\| X\right\| _{2}=\left( \sum_{ij}\left| X_{ij}\right|
^{2}\right) ^{1/2}.$ This gives us a Riemannian metric $ds$ on $\mathcal{U}%
\left( N\right) .$ The Riemannian metric on $\mathcal{SU}(N)$ can be defined
by restriction.

The (real or complex-valued) function $f\left( x\right) $ on a metric space $%
M$ is called Lipschitz with constant $L$ if for every two points $x,y\in M,$
it is true that $\left| f\left( x\right) -f\left( y\right) \right| \leq
Ld\left( x,y\right) $, where $d\left( x,y\right) $ is the shortest distance
between $x$ and $y.$

\begin{proposition}
\label{proposition_Blower}Let $g:\left( \mathcal{SU}\left( N\right) ,\left\|
ds\right\| _{2}\right) \rightarrow \mathbb{R}$ be an $L$-Lipschitz function
and let $\mathbb{E}g=0.$ Then \newline
(i) $\mathbb{E}\exp \left( tg\right) \leq \exp \left( ct^{2}L^{2}/N\right) $
for every $t\in \mathbb{R}$ and some numeric $c>0,$ and\newline
(ii) $P\left\{ \left| g\right| >\delta \right\} \leq \exp \left(
-c_{1}N\delta ^{2}/L^{2}\right) $ for every $\delta >0$ and some numeric $%
c_{1}>0.$
\end{proposition}

For the proof, see Theorems 3.8.3 and 3.9.2 in \cite{blower09} and Theorem
4.4.27 in \cite{anderson_guionnet_zeitouni10}.

In order to apply this result, we need to estimate the Lipschitz constants
for $m_{H}\left( z\right) $ and $f_{B}\left( z\right)$. If $M$ is a
Riemannian manifold and $f$ is a differentiable function on $M$, then it is
Lipschitz with constant $L$ provided that $\left| d_{X}f\left( x\right)
\right| \leq L$ for every $x\in M$ and every unit vector $X\in TM_{x}$. Here 
$d_{X}$ denotes the derivative in the direction of vector $X.$ We will apply
this general observation to the manifold $\mathcal{SU}(N)$.

Let $\widetilde{B}$ denote $UBU^{\ast },$ $B\left( X\right) =e^{iX}%
\widetilde{B}e^{-iX}$ and let 
\begin{equation*}
m_{H}\left( z,X\right) =\left( A+B\left( X\right) -z\right) ^{-1}
\end{equation*}%
We differentiate $m_{H}\left( z,X\right) $ with respect to $X$ (and evaluate
it at $X=0$) by using the chain rule.

\begin{eqnarray*}
\left| d_{X}m_{H}\left( z,X\right) \right| &=&\left| \sum_{x,y}\frac{%
\partial m_{H}\left( z\right) }{\partial \widetilde{B}_{xy}}%
d_{X}B_{xy}\left( X\right) \right| \\
&=&\left| \frac{1}{N}\sum_{x,y}\left( G^{2}\right) _{yx}\left[ X,\widetilde{B%
}\right] _{xy}\right| \\
&=&\left| \frac{1}{N}\sum_{x,y}\left( \left[ G^{2},\widetilde{B}\right]
\right) _{yx}X_{xy}\right| .
\end{eqnarray*}%
where we used the facts that $\partial m_{H}/\partial \left( \widetilde{B}%
_{xy}\right) =-N^{-1}\left( G^{2}\right) _{yx}$ and that $\left.
d_{X}B\left( X\right) \right| _{X=0}=\left[ X,\widetilde{B}\right] $. These
facts can be easily checked by a calculation. For the first one, see Lemma %
\ref{lemma_Gderivative_1} below.

If $\left\| X\right\| _{2}=1,$ then it follows that 
\begin{eqnarray*}
\left| d_{X}m\left( z,X\right) \right| &\leq &\frac{1}{N}\left\| \left[
G^{2},\widetilde{B}\right] \right\| _{2} \\
&\leq &\frac{1}{\sqrt{N}}\left\| \left[ G^{2},\widetilde{B}\right] \right\|
\\
&\leq &\frac{2\left\| B\right\| }{\sqrt{N}\eta ^{2}}
\end{eqnarray*}

Together with Proposition \ref{proposition_Blower}, this implies the first
claim of the lemma.

For the second claim, let $f_{B}\left( z,X\right) =B\left( X\right) \left(
A+B\left( X\right) -z\right) ^{-1}.$ Note that $\widetilde{B}\left( A+%
\widetilde{B}-z\right) ^{-1}=I-\left( A-z\right) \left( A+\widetilde{B}%
-z\right) ^{-1}.$ This allows us to calculate: 
\begin{equation*}
\frac{\partial }{\partial \widetilde{B}_{xy}}\left( f_{B}\left( z,X\right)
\right) =\frac{1}{N}\left( G\left( A-z\right) G\right) _{yx}.
\end{equation*}

Hence, 
\begin{eqnarray*}
\left| d_{X}f\left( z,X\right) \right| &=&\left| \sum_{x,y}\frac{\partial
f_{B}\left( z\right) }{\partial \widetilde{B}_{xy}}d_{X}B_{xy}\left(
X\right) \right| \\
&=&\left| \frac{1}{N}\sum_{x,y}\left( \left[ G\left( A-z\right) G,\widetilde{%
B}\right] \right) _{yx}X_{xy}\right| \\
&\leq &\frac{1}{N}\left\| \left[ G\left( A-z\right) G,\widetilde{B}\right]
\right\| _{2} \\
&\leq &\frac{1}{\sqrt{N}}\left\| \left[ G\left( A-z\right) G,\widetilde{B}%
\right] \right\| .
\end{eqnarray*}

Since $\left( A-z\right) G=I-$ $\widetilde{B}G,$ we can continue this as 
\begin{equation*}
\left| d_{X}f\left( z,X\right) \right| \leq \frac{2}{\sqrt{N}}\left( \frac{%
\left\| B\right\| }{\eta }+\frac{\left\| B\right\| ^{2}}{\eta ^{2}}\right) ,
\end{equation*}%
and the rest follows from Proposition \ref{proposition_Blower}. $\square $

Later, we will need the following consequence of Proposition \ref%
{theorem_m_large_deviations}.

\begin{corollary}
\label{corollary_sup_m_large_deviations}Let $I_{\eta }=[-2K+i\eta ,2K+i\eta
].$ Then for some positive $c$ and $c_{1}$ which may depend on $K$ and for
all $\delta >0$,%
\begin{equation*}
P\left\{ \sup_{z\in I_{\eta }}\left| m_{H}\left( z\right) -\mathbb{E}%
m_{H}\left( z\right) \right| >\delta \right\} \leq \exp \left( -\frac{%
c\delta ^{2}\eta ^{4}}{\left\| B\right\| ^{2}}N^{2}\right) ,
\end{equation*}%
provided that $N\geq c_{1}\left( \sqrt{-\log \left( \eta ^{2}\delta \right) }%
\right) /\left( \eta ^{2}\delta \right) .$
\end{corollary}

\textbf{Proof of Corollary:} Note that $\left| m_{H}^{\prime }\left(
z\right) \right| \leq \eta ^{-2}$ and $\left| \mathbb{E}m_{H}^{\prime
}\left( z\right) \right| \leq \eta ^{-2}$ and that it is enough to place $%
O\left( K/\eta ^{2}\delta \right) $ points on interval $I_{\eta }$ to create
an $\varepsilon $-net with $\varepsilon =\eta ^{2}\delta /4.$ If $\left|
m_{H}\left( z\right) -\mathbb{E}m_{H}\left( z\right) \right| \leq \delta /2$
at every point of the net, then $\left| m_{H}\left( z\right) -\mathbb{E}%
m_{H}\left( z\right) \right| \leq \delta $ for all $z\in I_{\eta }.$ Hence,
by Theorem \ref{theorem_m_large_deviations}, 
\begin{eqnarray*}
P\left\{ \sup_{z\in I_{\eta }}\left| m_{H}\left( z\right) -\mathbb{E}%
m_{H}\left( z\right) \right| >\delta \right\} &\leq &\frac{c^{\prime }K}{%
\eta ^{2}\delta }\exp \left( -\frac{c\delta ^{2}\eta ^{4}}{\left\| B\right\|
^{2}}N^{2}\right) \\
&=&\exp \left( -\frac{c\delta ^{2}\eta ^{4}}{\left\| B\right\| ^{2}}%
N^{2}+\log \left( \frac{c^{\prime }K}{\eta ^{2}\delta }\right) \right) \\
&\leq &\exp \left( -\frac{c^{\prime \prime }\delta ^{2}\eta ^{4}}{\left\|
B\right\| ^{2}}N^{2}\right) ,
\end{eqnarray*}%
if $N\geq c_{1}\left( \sqrt{-\log \left( \eta ^{2}\delta \right) }\right)
/\left( \eta ^{2}\delta \right) $ and $c_{1}$ is sufficiently large. $%
\square $

\section{An estimate on error terms in the Pastur-Vasilchuk system}

\label{Section_PV_System}

For the convenience of the reader, we re-derive here the Pastur-Vasilchuk
system. This is a system of equations for $\mathbb{E}m_{H}\left( z\right) $, 
$\mathbb{E}f_{A}\left( z\right) ,$ and $\mathbb{E}f_{B}\left( z\right) .$
When $N$ is large, this system is a perturbation of system (\ref%
{free_convolution_system}), and the main purpose of this section is to
estimate quantitatively the size of this perturbation. Later, we will show
that system (\ref{free_convolution_system}) is stable with respect to small
perturbations, and therefore for large $N$ the function $\mathbb{E}%
m_{H}\left( z\right) $ is close to the Stieltjes transform of $\mu
_{A_{N}}\boxplus \mu _{B_{N}}.$

We use notations 
\begin{equation*}
\Delta _{A}:=\left( m_{H}-\mathbb{E}m_{H}\right) G_{H}-G_{A}\left( f_{B}-%
\mathbb{E}f_{B}\right) G_{H}
\end{equation*}%
and 
\begin{equation}
R_{A}:=\frac{1}{\mathbb{E}m_{H}}\frac{1}{N}\mathrm{Tr}\left( \frac{1}{%
1+\left( \mathbb{E}f_{B}/\mathbb{E}m_{H}\right) G_{A}}\mathbb{E}\Delta
_{A}\right) ,  \label{definition_RA}
\end{equation}%
with similar definitions for $\Delta _{B}$ and $R_{B}$.

\begin{theorem}[Pastur-Vasilchuk]
\label{theorem_PV_system}The functions $\mathbb{E}m_{H}\left( z\right) ,$ $%
\mathbb{E}f_{A}\left( z\right) $ and $\mathbb{E}f_{B}\left( z\right) $
satisfy the following system of equations: 
\begin{eqnarray}
\mathbb{E}m_{H}\left( z\right) &=&m_{A}\left( z-\frac{\mathbb{E}f_{B}\left(
z\right) }{\mathbb{E}m_{H}\left( z\right) }\right) +R_{A}\left( z\right) ,
\label{system_PV} \\
\mathbb{E}m_{H}\left( z\right) &=&m_{B}\left( z-\frac{\mathbb{E}f_{A}\left(
z\right) }{\mathbb{E}m_{H}\left( z\right) }\right) +R_{B}\left( z\right) , 
\notag \\
z+\frac{1}{\mathbb{E}m_{H}\left( z\right) } &=&\frac{\mathbb{E}f_{A}\left(
z\right) +\mathbb{E}f_{B}\left( z\right) }{\mathbb{E}m_{H}\left( z\right) },
\notag
\end{eqnarray}%
where $R_{A}$ and $R_{B}$ are defined as in (\ref{definition_RA}).
\end{theorem}

The main technical tool in the proof of this theorem is the following
formula due to Pastur and Vasilchuk. Recall that $G_{H}$ is the resolvent of 
$H_{N}=A_{N}+U_{N}B_{N}U_{N}^{\ast }$ where $U_{N}$ is the Haar distributed
random unitary matrix.

\begin{proposition}
\label{proposition_main_identity}$\mathbb{E}\left( m_{H}G_{H}\right) =%
\mathbb{E}\left( m_{H}G_{A}-G_{A}f_{B}G_{H}\right) .$
\end{proposition}

This result immediately implies Theorem \ref{theorem_PV_system}. Indeed, the
identity in Proposition \ref{proposition_main_identity} can be written in
the following equivalent form.

\begin{eqnarray*}
(\mathbb{E}m_{H})\mathbb{E}G_{H} &=&(\mathbb{E}m_{H})G_{A}-(\mathbb{E}%
f_{B})G_{A}\mathbb{E}G_{H} \\
&+&\mathbb{E}[(m_{H}-\mathbb{E}m_{H})G_{H}]-G_{A}\mathbb{E}[(f_{B}-\mathbb{E}%
f_{B})G_{H}] \\
&=&(\mathbb{E}m_{H})G_{A}-(\mathbb{E}f_{B})G_{A}\mathbb{E}G_{H}+\mathbb{E}%
\Delta _{A}.
\end{eqnarray*}%
This expression can be further re-written (after we multiply it by $A_{N}-z$
and re-arrange terms) as 
\begin{equation*}
\mathbb{E}m_{H}\left( A_{N}-\left( z-\frac{\mathbb{E}f_{B}}{\mathbb{E}m_{H}}%
\right) \right) \mathbb{E}G_{H}=\mathbb{E}m_{H}+\left( A_{N}-z\right) 
\mathbb{E}\Delta _{A}.
\end{equation*}%
Let $z^{\prime }:=z-Ef_{B}/Em.$ Then for almost all values of $z,$ 
\begin{equation*}
\mathbb{E}m_{H}\mathbb{E}G_{H}=G_{A}\left( z^{\prime }\right) \mathbb{E}%
m_{H}+\left( A_{N}-z\right) G_{A}\left( z^{\prime }\right) \mathbb{E}\Delta
_{A}.
\end{equation*}%
Take the normalized trace and divide the resulting expression by $\mathbb{E}%
m_{H}.$ Then, we obtain 
\begin{eqnarray*}
\mathbb{E}m_{H}\left( z\right)  &=&m_{A}\left( z^{\prime }\right) +\frac{1}{%
\mathbb{E}m_{H}}\frac{1}{N}\mathrm{Tr}\left( \frac{1}{1+\left( \mathbb{E}%
f_{B}/\mathbb{E}m\right) G_{A}}\mathbb{E}\Delta _{A}\right) . \\
&=&m_{A}\left( z^{\prime }\right) +R_{A}.
\end{eqnarray*}%
The second equation of the system is obtained similarly and the third
equation is an identity.

\textbf{Proof of Prop. \ref{proposition_main_identity}: }It is useful to use
notation $\widetilde{B}=U_{N}B_{N}U_{N}^{\ast }$ and $B\left( X\right)
=e^{iX}\widetilde{B}e^{-iX}.$ Note that by using the resolvent identity $%
G_{H}\left( z\right) -G_{A}\left( z\right) =-G_{A}\left( z\right) \widetilde{%
B}G_{H}\left( z\right) ,$ we know that 
\begin{eqnarray*}
\mathbb{E}\left( m_{H}G_{H}\right) &=&\mathbb{E}\left( m_{H}G_{A}-m_{H}G_{A}%
\widetilde{B}G_{H}\right) \\
&=&G_{A}\mathbb{E}\left( m_{H}-m_{H}\widetilde{B}G_{H}\right) .
\end{eqnarray*}

Hence, it is enough to show that $\mathbb{E}\left( m_{H}\widetilde{B}%
G_{H}\right) =\mathbb{E}\left( f_{B}G_{H}\right) ,$

\begin{lemma}
\label{lemma_Gderivative_1}Let $A$ and $B$ be two arbitrary matrices and $%
G\left( z\right) =\left( A+B-z\right) ^{-1}.$ Then, $\left( \partial
G/\partial B_{xy}\right) _{uv}=-G_{ux}G_{yv}.$ In particular, 
\begin{equation*}
\left( \sum_{x,y}\left( \partial G/\partial B_{xy}\right) M_{xy}\right)
_{uv}=-\sum_{x,y}G_{ux}M_{xy}G_{yv}.
\end{equation*}
\end{lemma}

\textbf{Proof:} This is an immediate consequence of the resolvent identity $%
G_{X+Y}\left( z\right) -G_{X}\left( z\right) =-G_{X}\left( z\right)
YG_{X+Y}\left( z\right) $ applied to $X=A+B$ and $Y=tE^{xy}$, where $E^{xy}$
denote the matrix that have $1$ in the intersection of row $x$ and column $y$
and zeroes elsewhere. $\square$

\begin{lemma}
\label{lemma_Schwinger_Dyson_for_G} For every $u,v,a,b,$ it is true that 
\begin{equation*}
\mathbb{E}\left( \left( G_{H}\right) _{ua}\left( \widetilde{B}G_{H}\right)
_{bv}\right) =\mathbb{E}\left( \left( G_{H}\widetilde{B}\right) _{ua}\left(
G_{H}\right) _{bv}\right) .
\end{equation*}
\end{lemma}

\textbf{Proof}: Note that $d\left( \mathbb{E}\left[ \left( A+B\left(
X\right) -z\right) ^{-1}\right] \right) /dt=0$ for every Hermitian matrix $%
X, $ because the distribution of $B\left( X\right) =e^{-itX}\widetilde{B}%
e^{itX} $ is the same as the distribution of $\widetilde{B}.$ We can compute 
\begin{eqnarray*}
\left. \frac{d}{dt}\left[ \left( A+B\left( X\right) -z\right) ^{-1}\right]
\right| _{t=0} &=&\left. \sum_{x,y}\frac{\partial G_{H}}{\partial \widetilde{%
B}_{xy}}\frac{dB\left( X\right) _{xy}}{dt}\right| _{t=0} \\
&=&i\sum_{x,y}\frac{\partial G_{H}}{\partial \widetilde{B}_{xy}}\sum_{s}%
\left[ -X_{xs}\widetilde{B}_{sy}+\widetilde{B}_{xs}X_{sy}\right] .
\end{eqnarray*}

Let $E^{ab}$ denote an $N$-by-$N$ matrix that has zeros everywhere except at
the intersection of the $a$-th row and $b$-th column, where it has entry $1.$
If we set $X=E^{ab}+E^{ba}$ and use Lemma \ref{lemma_Gderivative_1}, then we
obtain 
\begin{eqnarray*}
-\mathbb{E}[(G_H) _{ua}(\widetilde{B}G_H)_{bv}+(G_H) _{ub}(\widetilde{B}G_H)
_{av}] \\
+\mathbb{E}[(G_H\widetilde{B})_{ua} (G_H)_{bv}+(G_H\widetilde{B}%
)_{ub}(G_H)_{av}] =0.
\end{eqnarray*}

If we set $X=i\left( E^{ab}-E^{ba}\right) ,$ then we obtain a similar
expression and adding them together, we get: 
\begin{equation*}
\mathbb{E}\left[ -\left( G_{H}\right) _{ua}\left( \widetilde{B}G_{H}\right)
_{bv}+\left( G_{H}\widetilde{B}\right) _{ua}\left( G_{H}\right) _{bv}\right]
=0.
\end{equation*}%
$\square$

If we take $u=a$ in the statement of Lemmas \ref{lemma_Schwinger_Dyson_for_G}%
$,$ then we get 
\begin{equation*}
\mathbb{E}\left( \left( G_{H}\right) _{aa}\left( \widetilde{B}G_{H}\right)
_{bv}\right) =\mathbb{E}\left( \left( G_{H}\widetilde{B}\right) _{aa}\left(
G_{H}\right) _{bv}\right) .
\end{equation*}%
By adding up these equalities over $a$ and dividing by $N,$ we obtain that $%
\mathbb{E}\left( m_{H}\widetilde{B}G_{H}\right) =\mathbb{E}\left(
f_{B}G_{H}\right) ,$ and Proposition \ref{proposition_main_identity} is
proved. $\square $

\bigskip

Now, we are going to estimate the error terms $R_{A}$ and $R_{B}.$ Let $\eta
_{0}\geq 0$ and $\kappa >0.$ For our purposes it is sufficient to make the
estimates in the region

\begin{equation*}
\Omega _{\eta _{0},\kappa }:=\left\{ z\in \mathbb{C}:\mathrm{Im}z\geq \eta
_{0},\mathrm{Im}z\geq \kappa \mathrm{Re}z\right\} .
\end{equation*}

\begin{proposition}
\label{proposition_estimate_error_term}Assume that $\max \left\{ \left\|
A\right\| ,\left\| B\right\| \right\} \leq K$ and let $\kappa >0.$ There
exists an $\eta _{0}$ $=cK$ such that for every $z=E+i\eta \in \Omega _{\eta
_{0},\kappa },$ it is true that 
\begin{equation*}
\left| R_{A}\right| \leq \frac{C}{N\eta ^{2}},
\end{equation*}%
where $C>0$ and depends only on $K$ and $\kappa .$
\end{proposition}

In order to prove this result, we will proceed in two steps. First, we will
estimate $\left\| \mathbb{E}\Delta _{A}\right\| .$ Then we estimate the
multipliers before $\mathbb{E}\Delta _{A}$ in the definition of $R_{A}.$

\begin{proposition}
\label{proposition_DeltaA_estimate}Let $z=E+i\eta .$ Assume that $\eta \geq
\eta _{0}$ and that $\max \left\{ \left\| A\right\| ,\left\| B\right\|
\right\} \leq K.$ Then 
\begin{equation*}
P\left\{ \left\| \Delta _{A}\left( z\right) \right\| \geq \varepsilon
\right\} \leq \exp \left[ -c\varepsilon ^{2}\eta ^{6}N^{2}\right] ,
\end{equation*}%
and $\left\| \mathbb{E}\Delta _{A}\left( z\right) \right\| \leq c/\left(
N\eta ^{3}\right) $ where constants depend only on $K$ and $\eta _{0}.$
\end{proposition}

\textbf{Proof of Proposition \ref{proposition_DeltaA_estimate}:}

\begin{lemma}
\label{lemma_estimate_1_part_Delta} Let $z=E+i\eta ,$ where $\eta >0.$ Then
for a numeric $c>0,$\newline
a)%
\begin{equation*}
P\left\{ \left\| \left( m_{H}\left( z\right) -\mathbb{E}m_{H}\left( z\right)
\right) G_{H}\right\| \geq \varepsilon \right\} \leq \exp \left[ -c\frac{%
\varepsilon ^{2}\eta ^{6}}{\left\| B\right\| ^{2}}N^{2}\right] ,
\end{equation*}%
and b) 
\begin{equation*}
P\left\{ \left\| G_{A}\left( f_{B}\left( z\right) -\mathbb{E}f_{B}\left(
z\right) \right) G_{H}\right\| \geq \varepsilon \right\} \leq \exp \left[ -c%
\frac{\varepsilon ^{2}\eta ^{8}}{\left\| B\right\| ^{4}}N^{2}/\left( 1+\frac{%
\eta }{\left\| B\right\| }\right) ^{2}\right] .
\end{equation*}
\end{lemma}

\textbf{Proof:} Note that if $X$ is a Hermitian matrix and $\eta >0$, then $%
\left\| \left( X-i\eta \right) ^{-1}\right\| \leq 1/\eta .$ By using this
fact and Proposition \ref{theorem_m_large_deviations}, we get 
\begin{equation*}
P\left\{ \left\| \left( m_{H}\left( z\right) -\mathbb{E}m_{H}\left( z\right)
\right) G_{H}\right\| \geq \delta /\eta \right\} \leq \exp \left[ -c\frac{%
\delta ^{2}\eta ^{4}}{\left\| B\right\| ^{2}}N^{2}\right] .
\end{equation*}%
Claim (a) of the lemma follows if we set $\varepsilon =\delta /\eta .$ Claim
(b) follows from Proposition \ref{theorem_m_large_deviations} in a similar
fashion. $\square $

The first claim of Proposition \ref{proposition_DeltaA_estimate} directly
follows from Lemma \ref{lemma_estimate_1_part_Delta}.

For the second claim, note that $\left\| \mathbb{E}\Delta _{A}\right\| \leq 
\mathbb{E}\left\| \Delta _{A}\right\| $ by the convexity of norm, and $%
\mathbb{E}\left\| \Delta _{A}\right\| $ can be estimated by using the first
claim of Proposition \ref{proposition_DeltaA_estimate} and the equality 
\begin{equation*}
\mathbb{E}X=\int_{0}^{\infty }\left( 1-\mathcal{F}_{X}\left( t\right)
\right) dt,
\end{equation*}%
valid for every positive random variable $X$ and its cumulative distribution
function $\mathcal{F}_{X}\left( t\right) .$ In our case, we obtain 
\begin{equation*}
\mathbb{E}\left\| \Delta _{A}\right\| \leq \int_{0}^{\infty }\exp \left[
-ct^{2}\eta ^{6}N^{2}\right] dt=\frac{c^{\prime }}{N\eta ^{3}}.
\end{equation*}%
$\square $

\begin{proposition}
\label{proposition_estimate_multiplier}Let $z=E+i\eta $ where $\eta >0.$
Assume that $\left\{ \left\| A\right\| ,\left\| B\right\| \right\} \leq K.$
Then, there exists such an $\eta _{0}=cK$ with numeric $c>0,$ that for every 
$\eta \geq \eta _{0},$ $\left\| \left( 1+\left( \mathbb{E}f_{B}/\mathbb{E}%
m_{H}\right) G_{A}\left( z\right) \right) ^{-1}\right\| \leq 2.$
\end{proposition}

The proof uses the following result.

\begin{lemma}
\label{lemma_power_expansions}Assume that $\left\{ \left\| A\right\|
,\left\| B\right\| \right\} \leq K.$ Then, for some numeric $c>0,$ the
functions $\mathbb{E}m\left( z\right) ,$ $\mathbb{E}f_{B}\left( z\right) ,$
and $\mathbb{E}m\left( z\right) /\mathbb{E}f_{B}\left( z\right) $ can be
represented by uniformly convergent series in $z^{-1}$ in the area $\left|
z\right| \geq cK,$ 
\begin{eqnarray*}
\mathbb{E}m\left( z\right) &=&-z^{-1}+\sum_{k=2}^{\infty }a_{k}\left[ m%
\right] z^{-k}, \\
\mathbb{E}f_{B}\left( z\right) &=&\sum_{k=1}^{\infty }a_{k}\left[ f_{B}%
\right] z^{-k}, \\
\frac{\mathbb{E}f_{B}\left( z\right) }{\mathbb{E}m\left( z\right) }
&=&\sum_{k=0}^{\infty }\beta _{k}z^{-k}.
\end{eqnarray*}
\end{lemma}

The proof of the first two equalities is by expansion of $\left(
A+B-z\right) ^{-1}$ and $B\left( A+B-z\right) ^{-1}$ in convergent series of 
$z^{-1}$ and estimating the coefficients in these series. This establishes
the uniform convergence in the area $\left| z\right| >cK$ and ensures that
it is possible to take expectation and trace of the series in a term-by-term
fashion. The third equality follows from the first two. $\square $

\textbf{Proof of Proposition \ref{proposition_estimate_multiplier}}: By the
previous lemma, $\mathbb{E}f_{B}/\mathbb{E}m_{H}$ is analytic in $z^{-1}$
and therefore bounded if $\left| z\right| >cK.$ Since $\left\| G_{A}\left(
z\right) \right\| \leq 1/\eta ,$ we can choose $\eta _{0}=cK$ with
sufficiently large $c,$ so that $\eta >\eta _{0}$ ensures that 
\begin{equation*}
\left\| \frac{\mathbb{E}f_{B}\left( z\right) }{\mathbb{E}m_{H}\left(
z\right) }G_{A}\left( z\right) \right\| <1/2,
\end{equation*}%
and 
\begin{equation*}
\left\| \left( 1+\frac{\mathbb{E}f_{B}\left( z\right) }{\mathbb{E}%
m_{H}\left( z\right) }G_{A}\left( z\right) \right) ^{-1}\right\| <2.
\end{equation*}%
$\square $

\textbf{Proof of Proposition \ref{proposition_estimate_error_term}:} For
every matrix $X,$ it is true that $\left| N^{-1}\mathrm{Tr}\left( X\right)
\right| \leq \left\| X\right\| $. Hence, by using Propositions \ref%
{proposition_DeltaA_estimate} and \ref{proposition_estimate_multiplier}, 
\begin{equation*}
\left| \frac{1}{N}\mathrm{Tr}\left( \frac{1}{1+\left( \mathbb{E}f_{B}/%
\mathbb{E}m_{H}\right) G_{A}}\mathbb{E}\Delta _{A}\right) \right| \leq
\left\| \frac{1}{1+\left( \mathbb{E}f_{B}/\mathbb{E}m_{H}\right) G_{A}}%
\right\| \left\| \mathbb{E}\Delta _{A}\right\| \leq \frac{c}{N\eta ^{3}},
\end{equation*}%
provided that $\eta >\eta _{0}=cK.$

By using the power expansion for $m(z)$, we find $m\left( z\right) ^{-1}\leq
2\left| z\right| \leq 2\eta \sqrt{1+\kappa ^{-2}}$ if $\left| z\right| >cK.$
It follows that for $z\in \Omega _{\eta _{0},\kappa },$%
\begin{equation*}
\left| R_{A}\right| \leq \frac{c}{N\eta ^{2}}\sqrt{1+\kappa ^{-2}}.
\end{equation*}%
$\square$

\section{Stability of the Pastur-Vasilchuk system}

\label{Section_stability_PV_system}

By results of \cite{pastur_vasilchuk00}, the solution of system (\ref%
{system_PV}) exists and unique in the upper half-plane $\mathbb{C}^{+}.$ We
are going to show that the solutions of systems (\ref{system_PV}) and (\ref%
{free_convolution_system}) are close to each other.

\begin{proposition}
\label{proposition_PV_stability}For all $z\in \Omega _{cK,\kappa },$ 
\begin{equation*}
\max \left\{ \left| \mathbb{E}m_{H}\left( z\right) -m_{\boxplus ,N}\left(
z\right) \right| \right\} \leq \frac{c^{\prime }}{N\eta },
\end{equation*}%
where $c$ and $c^{\prime }$ depends on $K$ and $\kappa $ only$.$
\end{proposition}

The idea of proof is to use the solution of the system (\ref%
{free_convolution_system}) as the starting point of the Newton-Kantorovich
algorithm (\cite{kantorovich48}) that computes the solution of system (\ref%
{system_PV}).

It is convenient to use a more uniform notation, so we write system (\ref%
{system_PV}) in a more compact form: 
\begin{eqnarray}
x_{1}-m_{A}\left( z-\frac{x_{3}}{x_{1}}\right) -R_{A} &=&0,
\label{system_PV_simplified} \\
x_{1}-m_{B}\left( z-\frac{x_{2}}{x_{1}}\right) -R_{B} &=&0,  \notag \\
zx_{1}-x_{2}-x_{3}+1 &=&0,  \notag
\end{eqnarray}%
The starting point of the algorithm is $x^{\boxplus }=\left( m_{\boxplus
,N},S_{A}m_{\boxplus ,N},S_{B}m_{\boxplus ,N}\right) ,$ where $m_{\boxplus
,N}\left( z\right) ,$ $S_{A}\left( z\right) ,$ and $S_{B}\left( z\right) $
are the solutions of (\ref{free_convolution_system}).The variable $z$ plays
the role of a parameter.

We assume that $R_{A}$ and $R_{B}$ are evaluated at the solution of (\ref%
{system_PV}) and fixed. Hence, in (\ref{system_PV_simplified}), $R_{A}$ and $%
R_{B}$ do not depend on $x$. The solution of (\ref{system_PV}) remains a
solution of this simplified system.

In a shorter form, system (\ref{system_PV_simplified}) can be written as \ 
\begin{equation}
P\left( x\right) =0.  \label{equation_Kantorovich}
\end{equation}%
Now, let us explain the Newton-Kantorovich method. Let (\ref%
{equation_Kantorovich}) be a general non-linear functional equation where $P$
is a non-linear operator that sends elements of a Banach space $X$ to a
Banach space $Y.$ Let $P$ be twice differentiable, and assume that the
operator $P^{\prime }\left( x\right) $ has an inverse $\left[ P^{\prime
}\left( x\right) \right] ^{-1}\in L\left( Y,X\right) $ where $L\left(
Y,X\right) $ denotes the space of bounded linear operators from $Y$ to $X.$
Then the Newton-Kantorovich method is given by the equation 
\begin{equation*}
x_{n+1}=x_{n}-\left[ P^{\prime }\left( x_{n}\right) \right] ^{-1}P\left(
x_{n}\right) .
\end{equation*}

The Kantorovich theorem (i) gives the sufficient conditions for the
convergence of this process, (ii) estimates the speed of convergence, and
(iii) estimates the distance of the solution $x^{\ast }$ from the initial
point $x_{0}.$ We give the statement of the theorem omitting the claim about
the speed of convergence, which is not important for us.

\begin{theorem}[Kantorovich]
\label{theorem_Kantorovich}Suppose that the following conditions hold:%
\newline
(1) for an initial approximation $x_{0},$ the operator $P^{\prime }\left(
x_{0}\right) \,$possesses an inverse operator $\Gamma _{0}=\left[ P^{\prime
}\left( x_{0}\right) \right] ^{-1}$ whose norm has the following estimate: $%
\left\| \Gamma _{0}\right\| \leq C_{0},$\newline
(2) $\left\| \Gamma _{0}P\left( x_{0}\right) \right\| \leq \delta _{0},$%
\newline
(3) the second derivative $P^{\prime \prime }\left( x\right) $ is bounded in
the domain determined by inequality (\ref{inequality_Kantorovich}) below;
namely, $\left\| P^{\prime \prime }\left( x\right) \right\| \leq M,$\newline
(4) the constants $C_{0},\delta _{0},M$ satisfy the relation $%
h_{0}=C_{0}\delta _{0}M\leq 1/2.$\newline
Then equation (\ref{equation_Kantorovich}) has a solution $x^{\ast },$ which
lies in a neighborhood of $x_{0}$ determined by the inequality 
\begin{equation}
\left\| x-x_{0}\right\| \leq \frac{1-\sqrt{1-2h_{0}}}{h_{0}}\delta _{0},
\label{inequality_Kantorovich}
\end{equation}%
and the successive approximations $x_{n}$ of the Newton method converge to $%
x^{\ast }.$
\end{theorem}

\textbf{Proof of Proposition \ref{proposition_PV_stability}:} In order to
apply the Newton-Kantorovich method, let us calculate the derivative $%
P^{\prime }\left( x\right) $ for our system: 
\begin{equation*}
P^{\prime }\left( x\right) =\left( 
\begin{array}{ccc}
1-m_{A}^{\prime }\left( z-\frac{x_{3}}{x_{1}}\right) \frac{x_{3}}{x_{1}^{2}}
& 0 & m_{A}^{\prime }\left( z-\frac{x_{3}}{x_{1}}\right) \frac{1}{x_{1}} \\ 
1-m_{B}^{\prime }\left( z-\frac{x_{2}}{x_{1}}\right) \frac{x_{2}}{x_{1}^{2}}
& m_{B}^{\prime }\left( z-\frac{x_{2}}{x_{1}}\right) \frac{1}{x_{1}} & 0 \\ 
z & -1 & -1%
\end{array}%
\right) .
\end{equation*}%
Then, the determinant is 
\begin{equation*}
\det \left( P^{\prime }\right) =-\frac{m_{A}^{\prime }+m_{B}^{\prime }}{x_{1}%
}+\frac{m_{A}^{\prime }m_{B}^{\prime }}{x_{1}^{3}}\left(
-zx_{1}+x_{2}+x_{3}\right) ,
\end{equation*}%
where $m_{A}^{\prime }$ and $m_{B}^{\prime }$ are short notations for $%
m_{A}^{\prime }\left( z-\frac{x_{3}}{x_{1}}\right) $ and $m_{B}^{\prime
}\left( z-\frac{x_{2}}{x_{1}}\right) ,$ respectively.

The power expansions from Lemma \ref{lemma_power_expansions} and the
definitions of $m_{A}$ and $m_{B}$ imply that $x_{1}\sim -z^{-1},$ $%
x_{2}\sim \alpha _{0}z^{-1},$ $x_{3}\sim \beta _{0}z^{-1},$ $m_{A}^{\prime
}\sim z^{-2},$ and $m_{B}^{\prime }\sim z^{-2}$ for $z\rightarrow \infty .$
Hence 
\begin{equation*}
\det \left( P^{\prime }\right) =\frac{1}{z}+O\left( 1\right) ,
\end{equation*}%
in the area $\left| z\right| >cK,$ where the constant in $O\left( 1\right) $
depends only on $K.$

(The proof that we gave for Lemma \ref{lemma_power_expansions} holds only
for $x_{1}=m_{H_{N}}\left( z\right) ,$ $x_{2}=f_{A_{N}}\left( z\right) ,$
and $x_{3}=f_{B_{N}}\left( Z\right) $. However, by using results from free
probability, these power expansions can be established in the case when $%
x_{1},$ $x_{2}$ and $x_{3}$ are defined as $m_{\boxplus ,N},$ $%
S_{A}m_{\boxplus ,N}$ and $S_{B}m_{\boxplus ,N},$ respectively.)

Now, it is easy to calculate the inverse of the derivative and find that 
\begin{equation}
\Gamma _{0}=\left[ P^{\prime }\left( x^{\boxplus }\right) \right]
^{-1}=z\left( 
\begin{array}{ccc}
0 & 0 & 0 \\ 
1 & 0 & 0 \\ 
0 & 1 & 0%
\end{array}%
\right) +O\left( 1\right) .  \label{formula_Gamma0}
\end{equation}

Hence 
\begin{eqnarray*}
\left\| \Gamma _{0}\right\| &=&\left| z\right| +O\left( 1\right) \\
&\leq &2\left| z\right| ,
\end{eqnarray*}%
if $z\in \Omega _{cK,\kappa }$ and $c$ is sufficiently large.

By using formula (\ref{formula_Gamma0}), we calculate for $z\in \Omega
_{cK,\kappa }$: 
\begin{eqnarray*}
\left\| \Gamma _{0}P\left( x^{\boxplus }\right) \right\| &\leq &\left|
z\right| \left( \left| R_{A}\right| +\left| R_{B}\right| \right) +O\left(
\left| R_{A}\right| +\left| R_{B}\right| \right) \\
&\leq &c\eta \left( \left| R_{A}\right| +\left| R_{B}\right| \right) \leq 
\frac{c^{\prime }}{N\eta },
\end{eqnarray*}%
where $c^{\prime }$ depends only on $K$ and $\kappa $ by Proposition \ref%
{proposition_estimate_error_term}.

The next step is to estimate $\left\| P^{\prime \prime }\left( x\right)
\right\| .$ Assume that $\left\| x-x^{\boxplus }\right\| \leq \frac{1}{2}%
\left| z\right| ^{-1}.$ (Later we will show that for large $N$ this disc
contains the disc given by (\ref{inequality_Kantorovich}).) By direct
computation of the second derivatives, it is easy to check that if $c$ is
sufficiently large and $z\in \Omega _{cK,\kappa },$ all second derivatives
of $P\left( x\right) $ are bounded by a constant, which can depend on $K$
only. Hence, $\left\| P^{\prime \prime }\left( x\right) \right\| \leq M,$
where $M$ depends on $K$ only.

Now we can apply Theorem \ref{theorem_Kantorovich} with $C_{0}=2\left|
z\right| ,$ $\delta _{0}=c^{\prime }/N\eta ,$ $M$ as in the previous
paragraph, and $h_{0}=C_{0}\delta _{0}M.$ For all sufficiently large $N,$ $%
h_{0}\leq 1/2$ and disc (\ref{inequality_Kantorovich}) is inside the disc $%
\left\| x-x^{\boxplus }\right\| \leq \frac{1}{2}\left| z\right| ^{-1}$ so
that the estimate for the second derivative holds.

Hence by Theorem \ref{theorem_Kantorovich}, if $z\in \Omega _{cK,\kappa },$
then the Newton algorithm which starts at $x^{\boxplus }$ will converge to a
solution of $P\left( x\right) =0$ and this solution satisfies inequality $%
\left\| x-x^{\boxplus }\right\| \leq 2\delta _{0}=c/\left( N\eta \right) .$
This completes the proof of Proposition \ref{proposition_PV_stability}. $%
\square$

\section{Hadamard's three circle theorem}

\label{section_hardy_theorem}

So far, we established the behavior of the difference $\left| \mathbb{E}%
m_{H_{N}}\left( z\right) -m_{\boxplus ,N}\left( z\right) \right| $ only for
the points where $\mathrm{Im}z\geq \eta _{0}.$ Here we prove a result about
its behavior for small $\mathrm{Im}z.$

\begin{proposition}
\label{proposition_bound_supremum_small_eta}Let $I_{\eta _{N}}$ be a
straight line segment between points $-2K+i\eta _{N}$ and $2K+i\eta _{N},$
where $\eta _{N}\geq c_{1}/\sqrt{\log N},$ and $c_{1}$ is a positive
constant that can depend on $K$. Then, 
\begin{equation*}
\sup_{z\in I_{\eta _{N}}}\left| \mathbb{E}m_{H}\left( z\right) -m_{\boxplus
,N}\left( z\right) \right| \leq \exp \left( -c\sqrt{\log N}\right) ,
\end{equation*}%
where $c$ depends only on $K.$
\end{proposition}

\begin{corollary}
\label{proposition_large_deviations_supm}Let $\eta _{N}=c_{1}\left( \log
N\right) ^{-\alpha },$ where $0<\alpha \leq 1/2$ and $I_{\eta _{N}}$ be a
straight line segment between points $-2K+i\eta _{N}$ and $2K+i\eta _{N}.$
Then, 
\begin{equation*}
P\left\{ \sup_{z\in I_{\eta _{N}}}\left| m_{H}\left( z\right) -m_{\boxplus
,N}\left( z\right) \right| >\delta \right\} \leq \exp \left( -c_{2}\delta
^{2}N^{2}\left( \log N\right) ^{-4\alpha }\right) .
\end{equation*}%
Constants $c_{1}$ and $c_{2}$ depend only on $K.$
\end{corollary}

\textbf{Proof of Corollary \ref{proposition_large_deviations_supm}:} This
result follows from Corollary \ref{corollary_sup_m_large_deviations} and
Proposition \ref{proposition_bound_supremum_small_eta}, which estimate $%
\left| m_{H}-\mathbb{E}m_{H}\right| $ and $\left| \mathbb{E}%
m_{H}-m_{\boxplus ,N}\right| ,$ respectively, if we note that for
sufficiently large $N,$ $\left| \mathbb{E}m_{H}\left( z\right) -m_{\boxplus
,N}\left( z\right) \right| <\delta $ for all $z\in I_{\eta _{N}}$. $\square $

For the proof of Proposition \ref{proposition_bound_supremum_small_eta}, we
use the three circle theorem by Hadamard (\cite{hardy15} or \cite{riesz23}).

\begin{theorem}[Hadamard's three circle theorem]
\label{theorem_Hardy}Suppose that $f\left( z\right) $ is a function of a
complex variable $z,$ holomorphic for $\left| z\right| <1,$ and let $M\left(
r\right) =\sup_{\theta }f\left( re^{i\theta }\right) $ for $r<1.$ Then $%
M\left( r\right) $ possesses the following properties:\newline
(1) $M\left( r\right) $ is an increasing function of $r$; \newline
(2) $\log M\left( r\right) $ is a convex function of $\log r,$ so that 
\begin{equation*}
\log M(r)\leq \frac{\log \left( r_{2}/r\right) }{\log \left(
r_{2}/r_{1}\right) }\log M(r_{2})+\frac{\log \left( r/r_{1}\right) }{\log
\left( r_{2}/r_{1}\right) }\log M(r_{1})
\end{equation*}%
if 
\begin{equation*}
0<r_{1}\leq r\leq r_{2}<1.
\end{equation*}
\end{theorem}

We will need the following consequence of this theorem.

\begin{lemma}
\label{lemma_consequence_Hardy}Suppose $f\left( z\right) $ is holomorphic
for $\left| z\right| <1,$ and let $M\left( r\right) $ be defined as in
Theorem \ref{theorem_Hardy}. Suppose that $M\left( r\right) \leq c/\left(
1-r\right) $ for all $r<1$ and that $M\left( e^{-1}\right) \leq \delta ,$
where $0<\delta <\delta _{0}$ and $\delta _{0}$ depends only on $c.$ Let 
\begin{eqnarray*}
r\left( \delta \right) &=&\exp \left( -4\sqrt{c/\log \left( 1/\delta \right) 
}\right) , \\
\varepsilon \left( \delta \right) &=&\exp \left( -\sqrt{c\log \left(
1/\delta \right) }\right) .
\end{eqnarray*}%
Then 
\begin{equation*}
M\left( r\right) \leq \varepsilon \left( \delta \right)
\end{equation*}%
for all $r\leq r\left( \delta \right) $.
\end{lemma}

(Note that if $\delta \rightarrow 0,$ then $r\left( \delta \right)
\rightarrow 1$ and $\varepsilon \left( \delta \right) \rightarrow 0.$)

\begin{figure}[tbph]
\begin{center}
\includegraphics[width=8cm]{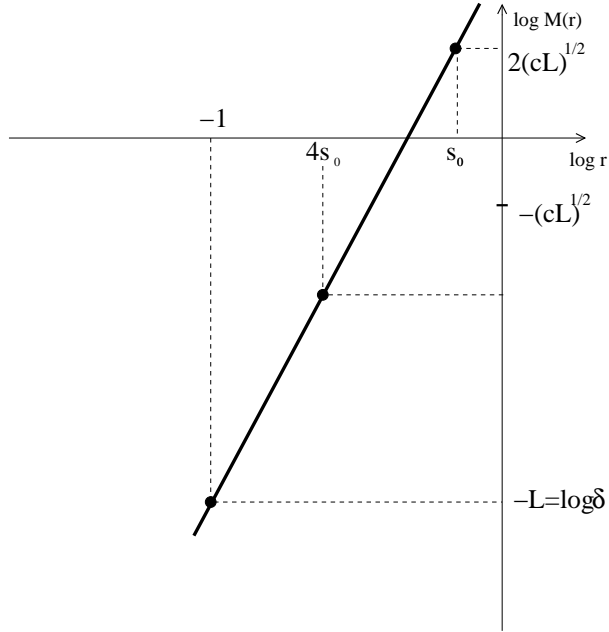}
\end{center}
\caption{Illustration to the proof of Lemma \ref{lemma_consequence_Hardy}}
\label{figure_AUBU_Hadamard_conseq}
\end{figure}

\textbf{Proof: } Let $L_{\delta }=\log \left( 1/\delta \right) $, $%
r_{0}=\exp \left( -\sqrt{c/L_{\delta }}\right) ,$ and $s_{0}=\log r_{0}=-%
\sqrt{c/L_{\delta }}.$ By assumption, 
\begin{equation*}
M\left( r_{0}\right) \leq \frac{c}{1-\exp \left( -\sqrt{c/L_{\delta }}%
\right) }\leq 2\sqrt{cL_{\delta }}
\end{equation*}%
for all sufficiently small $\delta .$ In the plane $\left( \log r,\log
M\right) $, the equation of the straight line that goes through points $%
\left( -1,-L_{\delta }\right) $ and $\left( \log r_{0},2\sqrt{cL_{\delta }}%
\right) $ is given by 
\begin{eqnarray*}
l\left( s\right) &=&\frac{2\sqrt{cL_{\delta }}+L_{\delta }}{-\sqrt{%
c/L_{\delta }}+1}\left( s+1\right) -L_{\delta } \\
&=&L_{\delta }\left[ \frac{2\sqrt{c}+\sqrt{L_{\delta }}}{-\sqrt{c}+\sqrt{%
L_{\delta }}}\left( s+1\right) -1\right] .
\end{eqnarray*}

By Hadamard's theorem, $\log M\left( e^{s}\right) \leq l\left( s\right) $
for all $s\in \left[ -1,s_{0}\right] .$ Let us set $\overline{s}=4s_{0}=-4%
\sqrt{c/L_{\delta }}.$ Then, 
\begin{equation*}
l\left( \overline{s}\right) =-\sqrt{cL_{\delta }}-\frac{9c\sqrt{L_{\delta }}%
}{-\sqrt{c}+\sqrt{L_{\delta }}}\leq -\sqrt{cL_{\delta }}
\end{equation*}%
if $L_{\delta }>c.$

Hence, 
\begin{equation*}
\log M\left( e^{s}\right) \leq -\sqrt{cL_{\delta }}
\end{equation*}%
if $s=\log r\leq \overline{s}=-4\sqrt{c/L_{\delta }}$ and $\delta \leq
\delta _{0}\left( c\right) .$ $\square $

Since we are interested in functions on the upper half-plane rather than on
the unit disc, we have to make a change of variables before we are able to
apply Hadamard's theorem. Consider the following map: 
\begin{equation*}
z=\frac{w-ia}{w+ia},
\end{equation*}%
where $a$ is a positive real number. This map sends the upper half-plane $%
\mathbb{C}^{+}=\left\{ w:\mathrm{Im}w\geq 0\right\} $ bijectively to the
unit disc $D\mathbb{=}\left\{ z:\left| z\right| \leq 1\right\} .$ In
particular, it sends point $ia$ to the center of the disc. The inverse
transformation is%
\begin{equation*}
w=ia\frac{1+z}{1-z}.
\end{equation*}%
Let $x\in \mathbb{R}$ and let $\xi =(x-ia)/\left( x+ia\right) \in \partial
D. $ Then 
\begin{equation*}
\frac{1}{x-w}=\frac{1}{2ai}\frac{\left( 1-\xi \right) \left( 1-z\right) }{%
\xi -z}.
\end{equation*}

Let 
\begin{equation*}
g\left( w\right) =\int_{-\infty }^{\infty }\frac{d\mu \left( x\right) }{x-w},
\end{equation*}%
where $\mathrm{Im}w>0$ and $\mu $ is the difference of two probability
measures. After the change of variable $w=w\left( z\right) ,$ this function
becomes a function of variable $z\in D.$ We will denote it as $f\left(
z\right) $. Then,%
\begin{equation}
f\left( z\right) =\frac{1}{2ai}\int_{\partial D}\frac{\left( 1-\xi \right)
\left( 1-z\right) }{\xi -z}d\nu \left( \xi \right) ,
\label{formula_Cauchy_in_circle}
\end{equation}%
where $\left| z\right| <1$ and $\nu $ is the forward image of $\mu ,$ hence
it is the difference of two probability measures on the unit circle $%
\partial D.$

Evidently, $f\left( z\right) $ is analytic for $\left| z\right| <1.$

\begin{lemma}
\label{lemma_parameter_bound}Let $f\left( z\right) $ be defined by formula (%
\ref{formula_Cauchy_in_circle}) with $\nu $ which is the difference of two
probability measures on $\partial D.$ Then, $M(r)\leq 4a^{-1}(1-r)^{-1}.$
\end{lemma}

\textbf{Proof:} Clearly, $\left| 1-\xi \right| \leq 2$, $\left| 1-z\right|
\leq 2,$ and $\left| \xi -z\right| \geq 1-\left| z\right| .$ It remains to
notice that the total variation of $\nu $ is bounded by 2, since it is a
difference of two probability measures. These facts imply that $\left|
f\left( z\right) \right| \leq 4a^{-1}\left( 1-\left| z\right| \right) ^{-1}.$
$\square$

\textbf{Proof of \ref{proposition_bound_supremum_small_eta}: }The map $w=ia%
\frac{1+z}{1-z}$ sends disc $B\left( 0,e^{-1}\right) $ to a disc $D_{1}\in 
\mathbb{C}^{+}$ that has the diameter 
\begin{equation*}
\left[ ia\frac{e-1}{e+1},ia\frac{e+1}{e-1}\right] .
\end{equation*}%
By an appropriate choice of $a,$ disc $D_{1}$ can be placed arbitrarily far
from the real axis, hence we can apply Proposition \ref%
{proposition_PV_stability} and write 
\begin{equation}
\sup_{w\in D_{1}}\left| \mathbb{E}m_{H}\left( w\right) -m_{\boxplus
,N}\left( w\right) \right| \leq \frac{c^{\prime }}{aN},
\label{estimate_on_D1}
\end{equation}%
where $c^{\prime }$ depends on $K.$

Next, define $\delta =c^{\prime }/\left( aN\right) $ and let $r\left( \delta
\right) =\exp \left( -8a^{-1}/\sqrt{\log \left( 1/\delta \right) }\right) $
as in Lemma \ref{lemma_consequence_Hardy} with parameter $c=4a^{-1}$. The
map $w=ia\frac{1+z}{1-z}$ sends disc $B\left( 0,r\left( \delta \right)
\right) $ to disc $D_{2}\in \mathbb{C}^{+}$ with the diameter 
\begin{equation*}
ia\left[ \frac{1-r\left( \delta \right) }{1+r\left( \delta \right) },\frac{%
1+r\left( \delta \right) }{1-r\left( \delta \right) }\right] .
\end{equation*}%
Note that the radius of $D_{2}$ approaches infinity as $\delta \downarrow 0,$
and that 
\begin{equation*}
ia\frac{1-r\left( \delta \right) }{1+r\left( \delta \right) }\sim 4i\sqrt{%
\frac{a}{\log \left( 1/\delta \right) }}=4i\sqrt{\frac{a}{\log \left(
aN/c^{\prime }\right) }}.
\end{equation*}%
It follows that there exists a $c_{1}>0$ such that for $\eta _{N}=c_{1}/%
\sqrt{\log N}$ all the points of the segment $I_{\eta _{N}}$ are located
inside the disc $D_{2}.$

Hence, Lemma \ref{lemma_consequence_Hardy} and estimate (\ref{estimate_on_D1}%
) imply that 
\begin{eqnarray}
\sup_{w\in I_{\eta _{N}}}\left| \mathbb{E}m_{H}\left( w\right) -m_{\boxplus
,N}\left( w\right) \right| &\leq &\exp \left( -2\sqrt{a^{-1}\log \left(
aN/c^{\prime }\right) }\right)  \label{estimate_on_I} \\
&\leq &\exp \left( -c_{2}\sqrt{\log N}\right) .
\end{eqnarray}

$\square$

\section{Proof of Theorems \ref{theorem_main} and \ref{theorem_local_law}}

\label{section_proof_main_theorem}

We use the following result due to Bai (see Theorems 2.1, 2.2, and Corollary
2.3 in \cite{bai93}). We formulate it in the form suitable for our
application

\begin{theorem}[Bai]
\label{theorem_bai}Let $K=\max \left\{ \left\| A_{N}.B_{N}\right\| \right\}
. $ Then, 
\begin{eqnarray}
\sup_{x}\left| \mathcal{F}_{H_{N}}\left( x\right) -\mathcal{F}_{\boxplus
,N}\left( x\right) \right| &\leq &c_{1}[\int_{-c_{2}K}^{c_{2}K}\left|
m_{H}\left( E+i\eta \right) -m_{\boxplus ,N}\left( E+i\eta \right) \right| dE
\notag \\
&&+\frac{1}{\eta }\sup_{E}\int_{\left| x\right| \leq 4\eta }\left| \mathcal{F%
}_{\boxplus ,N}\left( E+x\right) -\mathcal{F}_{\boxplus ,N}\left( E\right)
\right| dx],  \label{bai_inequality}
\end{eqnarray}%
where $c_{1}$ and $c_{2}$ are numeric.
\end{theorem}

\textbf{Proof of Theorem \ref{theorem_main}:} By using Assumption $A1,$ we
can estimate%
\begin{equation*}
\left| \mathcal{F}_{\boxplus ,N}\left( E+x\right) -\mathcal{F}_{\boxplus
,N}\left( E\right) \right| \leq T_{N}\left| x\right| ,
\end{equation*}
and therefore the second term on the right-hand side of (\ref{bai_inequality}%
) is bounded by $16T_{N}\eta .$

Let us set $\eta _{N}=c_{1}\left( \log N\right) ^{-\varepsilon /4},$ where $%
0<\varepsilon \leq 2.$ By Proposition \ref%
{proposition_bound_supremum_small_eta}, we can make 
\begin{equation*}
\sup_{z\in I_{\eta _{N}}}\left| \mathbb{E}m_{H}\left( z\right) -m_{\boxplus
,N}\left( z\right) \right| \leq \delta /3,
\end{equation*}%
provided that $N>\left( 3/\delta \right) ^{c\log \left( 3/\delta \right) }.$
We can also make $16T_{N}\eta _{N}\leq \delta /3$ by choosing $N\geq \exp
\left( \left( c/\delta \right) ^{4/\varepsilon }\right) .$

Then, we can use Bai's theorem and Corollary \ref%
{proposition_large_deviations_supm}, and find that for all sufficiently
large $N$ 
\begin{eqnarray*}
P\left\{ \sup_{x}\left| \mathcal{F}_{H_{N}}-\mathcal{F}_{\boxplus ,N}\right|
>\delta \right\} &\leq &P\left\{ \sup_{z\in I_{\eta _{N}}}\left| m_{H}\left(
z\right) -m_{\boxplus ,N}\left( z\right) \right| \geq c\delta \right\} \\
&\leq &P\left\{ \sup_{z\in I_{\eta _{N}}}\left| m_{H}\left( z\right) -%
\mathbb{E}m_{H}\left( z\right) \right| \geq c_{1}\delta \right\} \\
&\leq &\exp \left( -c_{2}\delta ^{2}N^{2}\left( \log N\right) ^{-\varepsilon
}\right) ,
\end{eqnarray*}%
where to make sure that the last inequality holds, it is enough to take 
\begin{equation*}
N\geq c_{1}\left( \sqrt{\log \left( 1/\left( \eta ^{2}\delta \right) \right) 
}\right) /\left( \eta ^{2}\delta \right) .
\end{equation*}%
For small $\delta ,$ the most binding inequality on $N$ is $N\geq \exp
\left( \left( c/\delta \right) ^{4/\varepsilon }\right) .$ $\square $

\bigskip

By using Theorem \ref{theorem_main}, we can derive the following corollary
and prove Theorem \ref{theorem_local_law}. Recall that $\mathcal{N}_{I}$
denotes the number of eigenvalues of $H$ in the interval $I.$

\begin{corollary}
\label{corollary_inequality_intervals}Suppose the assumptions of Theorem \ref%
{theorem_main} hold, and assume in addition that $\eta \geq c/(\varepsilon 
\sqrt{\log N})$. Then the following inequality holds: 
\begin{equation*}
P\left\{ \sup_{I,|I|=\eta }|\frac{\mathcal{N}_{I}}{N|I|}-\frac{\mu
_{\boxplus ,N}(I)}{|I|}|\geq \varepsilon \right\} \leq \exp \left(
-c\varepsilon ^{2}\frac{(\eta N)^{2}}{(\log N)^{2}}\right) ,
\end{equation*}%
where $c>0$ depends only on $K$ and $T$.
\end{corollary}

\textbf{Proof:} Let $I=\left( a,b\right] .$ Then $\mathcal{N}_{I}/N=\mathcal{%
F}_{H_{N}}\left( b\right) -\mathcal{F}_{H_{N}}\left( a\right) $ and $\mu
_{\boxplus ,N}\left( I\right) =\mathcal{F}_{\boxplus ,N}\left( b\right) -%
\mathcal{F}_{\boxplus ,N}\left( a\right) ,$ and therefore 
\begin{eqnarray*}
&&P\left\{ \sup_{I,|I|=\eta }\left| \frac{\mathcal{N}_{I}}{N|I|}-\frac{\mu
_{\boxplus ,N}(I)}{|I|}\right| \geq \varepsilon \right\} \\
&=&P\left\{ \sup_{a,b:b-a=\eta }\left| \mathcal{F}_{H_{N}}\left( b\right) -%
\mathcal{F}_{\boxplus ,N}\left( b\right) -\left( \mathcal{F}_{H_{N}}\left(
a\right) -\mathcal{F}_{\boxplus ,N}\left( a\right) \right) \right| \geq
\varepsilon \eta \right\} ,
\end{eqnarray*}%
and the corollary is the direct consequence of Theorem \ref{theorem_main}.
The assumption about $\eta $ is needed to ensure that $N$ in Theorem \ref%
{theorem_main} is sufficiently large and is forced by assumptions of
Proposition \ref{proposition_bound_supremum_small_eta}. $\square$

\textbf{Proof of Theorem \ref{theorem_local_law}: }Assumption $A1$ with
uniform $T$ ensures that $\mu _{\boxplus ,N}\left( I\right) /\left| I\right| 
$ approaches $\varrho _{\boxplus ,N}\left( E\right) $ when $I=\left( E-\eta
,E+\eta \right] $ and $\eta \rightarrow 0.$ Moreover, the convergence is
uniform in $E.$ Hence the conclusion of the theorem is implied by Corollary %
\ref{corollary_inequality_intervals}. $\square$

\section{Concluding remarks}

\label{section_conclusion}

We have shown that the probability of a large deviation of the empirical
c.d.f. of eigenvalues of $A_{N}+U_{N}B_{N}U_{N}^{\ast }\,$\ from the c.d.f.
of $\mu _{A_{N}}\boxplus \mu _{B_{N}}$ is bounded by $\exp \left( -c\delta
^{2}N^{2}/\log ^{\varepsilon }N\right) .$ The same results holds for the
ensemble in which $U_{N}$ denotes a Haar-distributed real orthogonal matrix.
In this case Lemma \ref{lemma_Schwinger_Dyson_for_G} does not hold as stated
and should be corrected. After this correction the identity in Proposition %
\ref{proposition_main_identity} becomes: 
\begin{equation*}
\mathbb{E}\left( m_{H}G_{H}\right) =\mathbb{E}\left(
m_{H}G_{A}-G_{A}f_{B}G_{H}\right) -\frac{1}{N}G_{A}\mathbb{E}\left( \left[
\left( G_{H}\right) ^{T},B\right] G_{H}\right) .
\end{equation*}%
Hence, we need to re-define $\Delta _{A}$ by adding an additional term 
\begin{equation*}
-N^{-1}G_{A}\left[ \left( G_{H}\right) ^{T},B\right] G_{H}.
\end{equation*}%
The norm of this term is bounded by $c/(N\eta ^{3}),$ therefore the estimate 
$\left\| E\Delta _{A}\right\| \leq c/\left( N\eta ^{3}\right) $ from
Proposition \ref{proposition_DeltaA_estimate} remains valid and further
analysis can be carried through without changes.

It would be interesting to investigate whether the empirical measure of
eigenvalues satisfies the large deviation principle. At the very least, it
should be expected that the limit 
\begin{equation*}
\lim_{N\rightarrow \infty }-\frac{1}{N^{2}}\log P\left\{ \left| \mathcal{F}%
_{H_{N}}\left( x\right) -\mathbb{E}\mathcal{F}_{H_{N}}\left( x\right)
\right| >\delta \right\}
\end{equation*}%
exists and is positive. It is also likely that the large deviation principle
holds at the level of measures. For classical Gaussian ensembles the large
deviation rate is closely related to the free entropy of a probability
measure: 
\begin{equation*}
\Sigma \left( \mu \right) =\int \log \left[ x-y\right] d\mu \left( x\right)
d\mu \left( y\right) .
\end{equation*}%
For more general large matrices with Gaussian entries, the large deviation
rates were obtained in the work of Guionnet. It is not clear if there are
similar formulas for the large deviation rate in the case of sums of random
matrices.

The second contribution of this paper is a local law for eigenvalues. It was
shown that the local law holds on the scale $(\log N)^{-1/2}$. It would be
interesting to extend this law to smaller scales. In the case when the
eigenvalue distributions of matrices $A_{N}$ and $B_{N}$ converge to
limiting distributions $\mu _{A}$ and $\mu _{B}$ with the free convolution $%
\mu _{A}\boxplus \mu _{B}$, the author expects that the local law holds on
the scale $N^{-1+\varepsilon }$ at all points where the density of the free
convolution exists. (A trivial cases when $\mu _{A}$ or $\mu _{B}$ are
concentrated on a single point should of course be ruled out.)

Currently, the limit laws on this scale are known for the Gaussian symmetric
and sample covariance matrices, where they are implied by the explicit
description of the limiting eigenvalue process on the scale $N^{-1}$. They
have also been established in \cite{erdos_schlein_yau_yin09} for the Wigner
and sample covariance random matrices. In this case, the local laws have
been used as the first step in the proof of the universality conjecture for
this class of random matrices.

Another area of possible further research is to understand better the local
structure of the eigenvalues, in particular, the point process of
eigenvalues and compare it to the structure of eigenvalues in classical
ensembles of random matrices. One would expect that the point process of
eigenvalues converges to a universal limit.

\bibliographystyle{plain}
\bibliography{comtest}

\begin{thebibliography}{10}

\bibitem{anderson_guionnet_zeitouni10}
Greg~W. Anderson, Alice Guionnet, and Ofer Zeitouni.
\newblock {\em An Introduction to Random Matrices}, volume 118 of {\em
  Cambridge studies in advanced mathematics}.
\newblock Cambridge University Press, 2009.

\bibitem{bai93}
Z.~D. Bai.
\newblock Convergence rate of expected spectral distributions of large random
  matrices. part \mbox{I}. \mbox{W}igner matrices.
\newblock {\em Annals of Probability}, 21:625--648, 1993.

\bibitem{belinschi08}
Serban~Teodor Belinschi.
\newblock The lebesgue decomposition of the free additive convolution of two
  probability distributions.
\newblock {\em Probability Theory and Related Fields}, 142:125--150, 2008.

\bibitem{ben_arous_guionnet97}
Gerard Be\mbox{n Arous} and Alice Guionnet.
\newblock Large deviations for \mbox{W}igner$\prime$s law and
  \mbox{V}oiculescu$\prime$s non-commutative entropy.
\newblock {\em Probability Theory and Related Fields}, 108:517--542, 1997.

\bibitem{bercovici_voiculescu98}
H.~Bercovici and D.~Voiculescu.
\newblock Regularity questions for free convolutions.
\newblock In H.~Bercovici and C.~Foias, editors, {\em Nonselfadjoint Operator
  Algebras, Operator Theory and Related Topics}, volume 104 of {\em Operator
  Theory Advances and Applications}, pages 37--47. Birkhauser: Basel, Boston,
  Berlin, 1998.

\bibitem{biane98b}
Philippe Biane.
\newblock Processes with free increments.
\newblock {\em Mathematische Zeitschrift}, 227:143--174, 1998.

\bibitem{blower09}
Gordon Blower.
\newblock {\em Random Matrices: High Dimensional Phenomena}, volume 367 of {\em
  London Mathematical Society Lecture Note Series}.
\newblock Cambridge University Press, 2009.

\bibitem{chatterjee07}
Sourav Chatterjee.
\newblock Concentration of \mbox{H}aar measures, with an application to random
  matrices.
\newblock {\em Journal of Functional Analysis}, 245:379--389, 2007.

\bibitem{erdos_schlein_yau_yin09}
L.~Erdos, B.~Schlein, H.-T. Yau, and J.~Yin.
\newblock The local relaxation flow approach to universality of the local
  statistics for random matrices.
\newblock preprint arXiv:0911.3687, 2009.

\bibitem{erdos_schlein_yau09}
Laszlo Erdos, Benjamin Schlein, and Horng-Tzer Yau.
\newblock Semicircle law on short scales and delocalization of eigenvectors for
  \mbox{W}igner random matrices.
\newblock {\em Annals of Probability}, 37:815--852, 2009.

\bibitem{gromov_milman83}
M.~Gromov and V.~D. Milman.
\newblock A topological application of isoperimetric inequality.
\newblock {\em American Journal of Mathematics}, 105(4):843--854, 1983.

\bibitem{hardy15}
G.~H. Hardy.
\newblock The mean value of the modulus of an analytic function.
\newblock {\em Proceedings of the London Mathematical Society}, 14:269--277,
  1915.

\bibitem{horn62}
A.~Horn.
\newblock Eigenvalues of sums of \mbox{H}ermitian matrices.
\newblock {\em Pacific Journal of Mathematics}, 12:225--241, 1962.

\bibitem{kantorovich48}
L.~V. Kantorovich.
\newblock Functional analysis and applied mathematics.
\newblock {\em Uspekhi Matematicheskih Nauk}, 3(6):89--185, 1948.
\newblock English translation available in L. V. Kantorovich, Selected Works,
  vol. 2, 171-280, (1996), Gordon and Breach Science Publishers.

\bibitem{knutson_tao99}
A.~Knutson and T.~Tao.
\newblock The honeycomb model of \mbox{$GL_n(\mathbb{C})$} tensor products
  \mbox{I}: Proof of the saturation conjecture.
\newblock {\em Journal of American Mathematical Society}, 12:1055--1090, 1999.

\bibitem{pastur_vasilchuk00}
L.~Pastur and V.~Vasilchuk.
\newblock On the law of addition of random matrices.
\newblock {\em Communications in Mathematical Physics}, 214:249--286, 2000.

\bibitem{riesz23}
F.~Riesz.
\newblock Sur les valeurs moyennes du module des fonctions harmonique et des
  fonctions analytiques.
\newblock {\em Acta Litterarum ac Scientiarum}, 1:27--32, 1922/23.
\newblock available in vol. 1 of the collected papers by F. Riesz.

\bibitem{speicher93}
Roland Speicher.
\newblock Free convolution and the random sum of matrices.
\newblock {\em Publications of RIMS (Kyoto University)}, 29:731--744, 1993.

\bibitem{voiculescu_dykema_nica92}
D.~Voiculescu, K.~Dykema, and A.~Nica.
\newblock {\em Free Random Variables}.
\newblock A.M.S. Providence, RI, 1992.
\newblock CRM Monograph series, No.1.

\bibitem{voiculescu91}
Dan Voiculescu.
\newblock Limit laws for random matrices and free products.
\newblock {\em Inventiones mathematicae}, 104:201--220, 1991.

\bibitem{weyl12}
H.~Weyl.
\newblock Das asymptotische \mbox{V}erteilungsgesetz der \mbox{E}igenwerte
  lineare partieller \mbox{D}ifferentialgleichungen.
\newblock {\em Mathematische Annalen}, 71:441--479, 1912.

\end{thebibliography}

\end{document}